\documentclass[a4paper,reqno,10pt,oneside]{amsart}
\usepackage{amsfonts,amssymb,latexsym,xspace,epsfig,color}

\usepackage{tikz}
\usetikzlibrary{arrows,snakes,backgrounds}
\usepackage{graphpap}
\usepackage{psfrag}
\usepackage{graphicx}
\usepackage[all]{xy}

\definecolor{BlueViolet}{cmyk}{0.86,0.91,0,0.04}
\definecolor{Black}{cmyk}{0,0,0,1}

\newtheorem{teo}{Theorem}

\newtheorem{lem}[teo]{Lemma}
\newtheorem{prop}[teo]{Proposition}

\newtheorem{defn}[teo]{Definition}
\newtheorem{exa}[teo]{Example}

\usepackage[latin1]{inputenc}
\usepackage{subfigure}
\def\bbR{{\mathbb R}}
\def\bbZ{{\mathbb Z}}
\def\bbN{{\mathbb N}}
\def\bbE{{\mathbb E}}
\def\bbV{{\mathbb V}}
\def\bbP{{\mathbb P}}

\def\T+{{\mathbb T_d^+}}

\def\K{{\mathcal K}}
\title{Evolution of a modified binomial random graph\\ by agglomeration}

\author{Mihyun Kang, Angelica Pachón and Pablo M. Rodríguez}

\address{
\newline
Mihyun Kang
\newline 
Institute of Discrete Mathematics, Graz University of Technology
\newline
Steyrergasse 30, 8010 Graz, Austria 
\newline
e-mail: kang@math.tugraz.at
\newline
\newline
 Angelica Pachón
 \newline
 Dipartimento di Matematica, Universit\`a di Torino, Via Carlo Alberto 10 - 10123 Torino, Italia
 \newline
 e-mail: angelicayohana.pachonpinzon@unito.it
 \newline
 \newline
 Pablo Martín Rodríguez
 \newline 
 Instituto de Ci\^encias Matem\'aticas e de Computa\c{c}\~ao, Universidade de S\~ao Paulo
 \newline 
 Av. Traba\-lha\-dor S\~ao-Carlense 400 - Centro, CEP 13560-970, S\~ao Carlos, SP, Brazil
 \newline
 Current Address: Laboratoire de Probabilités et Modèles Aléatoires, Université Paris-Diderot
 \newline
 Bâtiment Sophie Germain, Avenue de France - 75013 Paris, France
 \newline 
 e-mail: pablor@icmc.usp.br
}

\subjclass[2010]{05C80, 60C05}
\keywords{Erd\H os-R\'enyi model, Random Graph, Inhomogeneous Random Graph, Connectedness, Phase Transition}

\begin{document}

%\maketitle

\begin{abstract}
In the classical Erd\H os-R\'enyi random graph $G(n,p)$ there are $n$ vertices and each of the possible edges is independently present with probability $p$. The random graph $G(n,p)$ is homogeneous in the sense that all vertices have the same characteristics. On the other hand, numerous real-world networks are inhomogeneous in this respect. Such an inhomogeneity of vertices may influence the connection probability between pairs of vertices. %In order to capture typical properties of real-world networks,  various random graph models have recently been introduced and analyzed.

The purpose of this paper is to propose a new inhomogeneous random graph model which is obtained in a constructive way from the Erd\H os-R\'enyi random graph $G(n,p)$. Given a configuration of $n$ vertices arranged in $N$ subsets of vertices (we call each subset a super-vertex), we define a random graph with $N$ super-vertices  by letting two super-vertices be connected if and only if there is at least one edge between them in $G(n,p)$. Our main result concerns the threshold for connectedness. We also analyze the  phase transition for the emergence of the giant component and the degree distribution. 

Even though our model begins with $G(n,p)$, it assumes the existence of some community structure encoded in the configuration. Furthermore, under certain conditions it exhibits a power law degree distribution. Both properties are important for real applications.
\end{abstract}

\maketitle

%\begin{keywords}
%Erd\H os-R\'enyi Model, Random Graph, Inhomogeneous Random Graph, Connectedness, Phase Transition
%\end{keywords}

%\begin{AMS}
%05C80, 60C05
%\end{AMS}

%%%%%%%%%%%%%%%%%%%%%%%%%%%%%%%%%%%%%%%%%%%%%%%%%%%%%%%%%%%%%%%%%%%%%%%%%%%%%%
%%%%%%%%%%%%% INTRODUCTION
%%%%%%%%%%%%%%%%%%%%%%%%%%%%%%%%%%%%%%%%%%%%%%%%%%%%%%%%%%%%%%%%%%%%%%%%%%%%%%

\pagestyle{myheadings}
\thispagestyle{plain}
\markboth{M. KANG ET AL.}{A MODIFIED BINOMIAL RANDOM GRAPH}

\section{Introduction}
\label{sec:in}
The subject of random graphs began in 1959-1960 with the papers ``On random graphs I'' and  ``On the evolution of random graphs'' by Erd\H os and R\'enyi \cite{ERI,ER}. Since then, many properties of the Erd\H os-R\'enyi random graph  have been analyzed in order to answer questions of mathematical and physical interest. 

The original model studied by  Erd\H os and R\'enyi is the uniform random graph $G(n,M)$,  which is a  graph chosen uniformly at random among all graphs with vertex set $[n]:=\{1,2,\ldots,n\}$ and exactly $M$ edges.  It is not difficult to see that $G(n,M)$ is closely related to the the  binomial random graph $G(n,p)$, which is a graph with vertex set $[n]$, in which each pair of vertices is connected by an edge with probability $p$,  independently of each other. The latter model  was introduced by Gilbert \cite{Gilbert} at about the same time. It is well known that the two random graph models $G(n, p)$ and $G(n, M)$ are essentially equivalent for the correct choice of $M$ and $p$. Due to this equivalence and the deep and notable results proved in \cite{ERI,ER}, the binomial random graph $G(n,p)$ is also known in literature as the Erd\H os-R\'enyi random graph. 

One of the most classical results on the  Erd\H os-R\'enyi random graph is the threshold for connectedness, which is closely related to the non-existence of isolated vertices. Clearly,  when there exists at least one isolated vertex the graph is disconnected, but the opposite implication is not  generally true. Remarkably, it was shown that when there are no isolated vertices, the random graph $G(n, p)$ is connected with high probability (in short $whp$) which means with probability tending to one as $n$ goes to $\infty$  (see \cite{BT, ER}).

Many other properties were also studied in \cite{ER}. One of the most striking results was the discovery of a drastic change in the size of the largest component when the number of edges passes through $n/2$.  This phenomenon is related to the phase transition in percolation, a model that is well studied in mathematical physics and  is  one of the main branches of contemporary probability, see \cite{BollobasRiordan,Grimmett}.   However, perhaps the first to talk about  this kind of phenomena were Flory and Stockmayer,  using a very different language   coming from polymer physics,  concerning gelation rather than percolation, and a gel rather than a giant component, see for instance \cite{PJFlory}.

During the last few decades, increasing interest in the field of random graphs has been devoted to find models that describe the complexity of real-world networks. It has recently been observed that many real-world networks are inhomogeneous, in the sense that nodes may be of different types and their connections may depend on types (see, for example, \cite{newman}). A general theoretical model of an inhomogeneous random graph is proposed in the seminal paper of Bollob\'as, Janson and Riordan \cite{BJR}, who considered a conditional independence between the edges, where the number of edges is linear in the number of vertices. This model includes as special cases many models previously studied in literature, for instance Durrett \cite{Durret} and Bollob\'as, Janson and Riordan \cite{BJR2}. In the inhomogeneous random graph model introduced in \cite{BJR}, it is shown that under a weak (convergence) assumption on the expected number of edges, many interesting properties can be determined, in particular the critical point of the phase transition and the size of the giant component. More recently, van der Hofstad \cite{Hofstad} analyzed the critical behavior of the largest component in inhomogeneous random graphs in the so-called rank-1 case, where weights are associated with the vertices of the graph, and edges are present between vertices with a probability that is approximately proportional to the product of the respective weights. 

In this paper we propose a new inhomogeneous random graph model that is obtained in a constructive way from the classical Erd\H os-R\'enyi model. By ``a constructive way'', we mean an explicit scheme for constructing the graph from a given realization of $G(n,p)$. Given a partition of the vertex set $[n]$ of $G(n,p)$, in which each partition class represents an agglomeration of nodes, we call each partition class a super-vertex and the partition a configuration of super-vertices. We define an inhomogeneous random graph model by letting two super-vertices be connected if and only if there is at least one edge between them in $G(n,p)$. Note that our model assumes the existence of a kind of community structure by the agglomeration of the nodes, which is encoded in the super-vertices. However, we are not assuming that the vertices inside each super-vertex should be all connected. In other words, each super-vertex is not necessarily a clique in $G(n,p)$ and can be any kind of subgraph of $G(n,p)$. Related random graph models are analyzed, for instance, by Janson and Spencer \cite{JansonSpencer} and by Seshadhri, Kolda and Pinar \cite{SKP}.  

The main contributions of this paper are fourfold. We determine  (i) the threshold for the connectedness of our inhomogeneous random graph model (Theorem \ref{connectedness}); (ii) the threshold for the existence of the giant component formed by super-vertices (Proposition \ref{corBollobas1}); (iii) the degree distribution of a super-vertex (Proposition \ref{corBollobas}). Finally we show that under certain conditions our model exhibits a {\em power law degree distribution} (Example \ref{powerlaw}), which is an important property for real applications.

In order to determine the asymptotic probability of our model  being connected, we analyze the distribution of the number of isolated super-vertices, using the second moment method (in Lemma \ref{zeroone-iso}), as well as Stein's method (in Lemma \ref{poisson}). As for the threshold for the existence of the giant component  and the degree distribution of  super-vertices, we show that our model can be viewed as a special case of the inhomogeneous random graph (IRG) studied in \cite{BJR}, by identifying a graphical sequence of kernels  and applying the corresponding results for IRG in \cite{BJR}.

The rest of the paper is organized as follows. In Section \ref{S: The model} we introduce our inhomogeneous random graph model, which is followed by our main results and related work. The proof of  the location of the threshold for connectedness  (Theorem \ref{connectedness}) is provided in Section \ref{proofcon}. The proofs for  the emergence of the giant component and the degree distribution are given in Section \ref{graphicalproof}.  Finally,  in Section \ref{sec:diss} we add some concluding remarks. 

%%%%%%%%%%%%%%%%%%%%%%%%%%%%%%%%%%%%%%%%%%%%%%%%%%%%%%%%%%%%%%%%%%%%%%%%%%%%%%
%%%%%%%%%%%%% THE MODEL AND RESULTS
%%%%%%%%%%%%%%%%%%%%%%%%%%%%%%%%%%%%%%%%%%%%%%%%%%%%%%%%%%%%%%%%%%%%%%%%%%%%%%

\bigskip
\section{Our model and main results}
\label{S: The model}

\subsection{The model} 
The motivation of our model comes from real-world networks that exhibit a {\em community structure}. In order to reflect a possible community structure  we shall define a random graph with a given number of  super-vertices (also known as agglomerates) of given sizes, in which we assume that the underlying graph follows the $G(n,p)$-law. 

More precisely speaking, for each $N\in\bbN:=\{1,2,\ldots\}$ we let $r\in \bbN\cup\{\infty\}$ be either a constant independent of $N$ or a function in $N$ such that $r=r(N)$ tends to a constant or $\infty$ as $N\rightarrow \infty$. We let $\K^r:=\{(k_1,\ldots,k_r) \in\bbN^r: \sum_{i=1}^r k_{i}=N\}$ and $p=p(N,\K^r)\in[0,1]$ be given. Note that $\K^r$ and $k_1,\ldots,k_r$  depend on $N$, but for the sake of simplicity, we suppress this dependence in our notation. We may use the notation $\K^{\infty}$  when there exists a super-vertex of size $r=r(N)$ satisfying  $r(N)\rightarrow \infty$ as $N\rightarrow \infty$. 

Given a partition of the vertex set $[n]$ of $G(n,p)$,  we call a partition class of size $i$ {\em a super-vertex of size $i$}. Note that the vertices in each super-vertex are {\em not necessarily connected} in $G(n,p)$. We define  $G(N,\K^r,p)$ to be a random graph with $N$ super-vertices with configuration $\K^r$, in which  for each $i,j=1,2,\ldots,r$ there are $k_i$ super-vertices of size $i$ and an edge between a pair of two distinct super-vertices of sizes $i$ and $j$ is present with probability 
 \begin{eqnarray}\label{prob}
p_{ij} 
 & := & 1-\left(1-p\right)^{ij},
\end{eqnarray}      
independently of each other. In words, $p_{ij} $ is the probability that there is at least one edge between the corresponding partition classes of the vertex set $[n]$ of $G(n,p)$, see Figure 1. Note that the number of super-vertices and the number of vertices   are given by $N=\sum_{i=1}^{{r}}k_{i}$ and $n=\sum_{i=1}^{r}ik_{i}$, respectively.
 
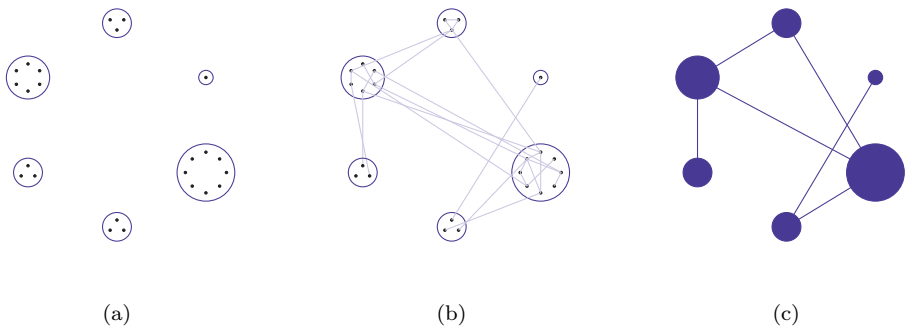
\begin{figure}[h]\label{fig}
\begin{center}
\begin{tabular}{ccc}

\begin{tikzpicture}[scale=0.9]

\draw [BlueViolet] (1.3,1.2) circle (3pt);
\filldraw [Black] (1.3,1.2) circle (0.5pt);

\draw [BlueViolet] (0,2) circle (6pt);
\filldraw [Black] (0.1,2.05) circle (0.5pt);
\filldraw [Black] (-0.1,2.05) circle (0.5pt);
\filldraw [Black] (0,1.9) circle (0.5pt);

\draw [BlueViolet] (-1.3,1.2) circle (9pt);
\filldraw [Black] (-1.3,1.4) circle (0.5pt);
\filldraw [Black] (-1.3,1) circle (0.5pt);
\filldraw [Black] (-1.47,1.3) circle (0.5pt);
\filldraw [Black] (-1.13,1.3) circle (0.5pt);
\filldraw [Black] (-1.47,1.1) circle (0.5pt);
\filldraw [Black] (-1.13,1.1) circle (0.5pt);

\draw [BlueViolet] (-1.3,-0.2) circle (6pt);
\filldraw [Black] (-1.3,-0.1) circle (0.5pt);
\filldraw [Black] (-1.2,-0.25) circle (0.5pt);
\filldraw [Black] (-1.4,-0.25) circle (0.5pt);

\draw [BlueViolet] (0,-1) circle (6pt);
\draw (0,-2) node[below,font=\footnotesize] {(a)};
\filldraw [Black] (0,-0.9) circle (0.5pt);
\filldraw [Black] (0.1,-1.05) circle (0.5pt);
\filldraw [Black] (-0.1,-1.05) circle (0.5pt);

\draw [BlueViolet] (1.3,-0.2) circle (12pt);
\filldraw [Black] (1.5,-0.4) circle (0.5pt);
\filldraw [Black] (1.1,0) circle (0.5pt);
\filldraw [Black] (1.5,0) circle (0.5pt);
\filldraw [Black] (1.1,-0.4) circle (0.5pt);
\filldraw [Black] (1.3,-0.5) circle (0.5pt);
\filldraw [Black] (1.3,0.1) circle (0.5pt);
\filldraw [Black] (1.6,-0.2) circle (0.5pt);
\filldraw [Black] (1,-0.2) circle (0.5pt);

\end{tikzpicture}

\hspace{0.8cm}&\hspace{0.8cm}
\begin{tikzpicture}[scale=0.9]

\draw [BlueViolet] (1.3,1.2) circle (3pt);
\filldraw [Black] (1.3,1.2) circle (0.5pt);

\draw [BlueViolet] (0,2) circle (6pt);
\filldraw [Black] (0.1,2.05) circle (0.5pt);
\filldraw [Black] (-0.1,2.05) circle (0.5pt);
\filldraw [Black] (0,1.9) circle (0.5pt);

\draw[BlueViolet!20] (0.1,2.05) -- (0,1.9);
\draw[BlueViolet!20] (0.1,2.05) -- (-0.1,2.05);

\draw [BlueViolet] (-1.3,1.2) circle (9pt);
\filldraw [Black] (-1.3,1.4) circle (0.5pt);
\filldraw [Black] (-1.3,1) circle (0.5pt);
\filldraw [Black] (-1.47,1.3) circle (0.5pt);
\filldraw [Black] (-1.13,1.3) circle (0.5pt);
\filldraw [Black] (-1.47,1.1) circle (0.5pt);
\filldraw [Black] (-1.13,1.1) circle (0.5pt);

\draw[BlueViolet!20] (-1.3,1.4) -- (-1.13,1.1);
\draw[BlueViolet!20] (-1.47,1.1) -- (-1.47,1.3);
\draw[BlueViolet!20]  (-1.3,1) -- (-1.13,1.3);
\draw[BlueViolet!20]  (-1.13,1.1) -- (0,1.9);
\draw[BlueViolet!20] (-1.47,1.1) -- (-1.2,-0.25);
\draw[BlueViolet!20]  (-1.3,1) -- (-1.3,-0.1);
%\draw[BlueViolet!20] (-1.3,1) -- (-1.2,-0.25);
\draw[BlueViolet!20]  (-1.47,1.3) -- (0,1.9);

\draw [BlueViolet] (-1.3,-0.2) circle (6pt);
\filldraw [Black] (-1.3,-0.1) circle (0.5pt);
\filldraw [Black] (-1.2,-0.25) circle (0.5pt);
\filldraw [Black] (-1.4,-0.25) circle (0.5pt);

\draw [BlueViolet] (0,-1) circle (6pt);
\draw (0,-2) node[below,font=\footnotesize] {(b)};
\filldraw [Black] (0,-0.9) circle (0.5pt);
\filldraw [Black] (0.1,-1.05) circle (0.5pt);
\filldraw [Black] (-0.1,-1.05) circle (0.5pt);

\draw[BlueViolet!20] (0,-0.9) -- (1.3,1.2);

\draw [BlueViolet] (1.3,-0.2) circle (12pt);
\filldraw [Black] (1.5,-0.4) circle (0.5pt);
\filldraw [Black] (1.1,0) circle (0.5pt);
\filldraw [Black] (1.5,0) circle (0.5pt);
\filldraw [Black] (1.1,-0.4) circle (0.5pt);
\filldraw [Black] (1.3,-0.5) circle (0.5pt);
\filldraw [Black] (1.3,0.1) circle (0.5pt);
\filldraw [Black] (1.6,-0.2) circle (0.5pt);
\filldraw [Black] (1,-0.2) circle (0.5pt);

\draw[BlueViolet!20] (1.3,0.1) -- (1.3,-0.5);
\draw[BlueViolet!20] (1.1,-0.4) -- (1,-0.2);
\draw[BlueViolet!20] (1.1,0) -- (1,-0.2);
\draw[BlueViolet!20] (1.5,-0.4) -- (1.6,-0.2);
\draw[BlueViolet!20] (1.5,0) -- (1.1,-0.4);
\draw[BlueViolet!20] (1.1,-0.4) -- (-1.47,1.3);
\draw[BlueViolet!20] (1.1,0) -- (-1.13,1.1);
\draw[BlueViolet!20]  (1.1,0) -- (0.1,-1.05);
\draw[BlueViolet!20] (1.3,0.1) -- (-0.1,2.05);
\draw[BlueViolet!20] (1.3,-0.5) -- (1.1,0);
\draw[BlueViolet!20] (1.3,0.1) -- (1.1,0);
\draw[BlueViolet!20]  (1.6,-0.2) -- (1.1,0);
\draw[BlueViolet!20]  (1.6,-0.2) -- (-1.13,1.3);
\draw[BlueViolet!20] (1.3,-0.5) -- (-0.1,-1.05);
\draw[BlueViolet!20] (1.3,0.1) -- (-1.3,1);

\end{tikzpicture}
\hspace{0.8cm}&\hspace{0.8cm}

\begin{tikzpicture}[scale=0.9]

\filldraw [BlueViolet] (1.3,1.2) circle (3pt);

\filldraw [BlueViolet] (0,2) circle (6pt);
\draw (0,-2) node[below,font=\footnotesize] {(c)};

\filldraw [BlueViolet] (-1.3,1.2) circle (9pt);

\filldraw [BlueViolet] (-1.3,-0.2) circle (6pt);

\filldraw [BlueViolet] (0,-1) circle (6pt);

\filldraw [BlueViolet] (1.3,-0.2) circle (12pt);

\draw[BlueViolet] (1.3,-0.2) -- (0,2);
\draw[BlueViolet] (1.3,-0.2) -- (-1.3,1.2);
\draw[BlueViolet] (1.3,-0.2) -- (0,-1);
\draw[BlueViolet] (1.3,1.2) -- (0,-1);
\draw[BlueViolet] (-1.3,-0.2) -- (-1.3,1.2);
\draw[BlueViolet] (-1.3,1.2) -- (0,2);

\end{tikzpicture}
\end{tabular}
\end{center}
\caption{Construction of $G(N,\K^r,p)$. (a) Begin with a fixed configuration $\K^r$ of $N$ super-vertices, that is, subsets of vertices (there are $n$ vertices in total). (b) Connect every pair of vertices independently with probability $p$,  in other words, take a realization of $G(n,p)$.
(c) A pair of super-vertices in $G(N,\K^r,p)$ is connected if and only if there is at least one edge between the corresponding subsets of vertices in $G(n,p)$.}
\end{figure}
 
Our goal is to study  properties of the random graph $G(N,\K^r,p)$ by considering different values of $p$ as a function of both $n$ and $N$. It seems difficult to obtain substantial results for $G(N,\K^r,p)$ without further restrictions. Throughout the paper we therefore assume that for each $i=1,2,\ldots,r$ the following limit exists 
\begin{equation}\label{convki}
\mu_i:=\lim_{N\to \infty}\frac{k_{i}}{N}
\end{equation}
and that $\mu_i>0 \text{ for some }i\in\{1,2,\dots,r\}$ (which means that for $N$ sufficiently large, there are linearly many super-vertices of size $i$ for some $i\in\{1,2,\dots,r\}$). 
 Furthermore, we define 
\begin{equation}\label{maxi_*} 
i_*=i_*(N):=\Big\{ 1\le i \le r \ \Big|\  k_{i}^{\frac{1}{i (n-i)}} = \max_{1\le j\le r} k_{j}^{\frac{1}{j (n-j)}}\Big\}.
\end{equation}
The  quantity $i_*$  will turn out to be crucial for the existence of isolated super-vertices and therefore for the threshold for connectedness.

Observe that {\em if} $r$ is either a constant independent of $N$ or a function in $N$ such that $r=r(N)$ tends to a constant as $N\to \infty$  (and so we may assume without loss of generality that $1\le r\le r^*$  for some constant $r^*\in \mathbb N$ for $N$ sufficiently large), then  $i_*$ is bounded above by a constant independent of $N$ and \eqref{convki} implies that $\lim_{N\to \infty}n/N$ exists and  $\lim_{N\to \infty}n/N=\sum_{i=1}^r i \mu_i$.

\bigskip
\subsection{Main results} 
The first main result of this paper concerns the exact location of the threshold for connectedness in $G(N,\K^r,p)$. The property of connectedness has not been addressed in the study of the general inhomogeneous random graph model (IRG) by Bollob\'as, Janson and Riordan \cite{BJR}. For a special class of IRG, the connectedness has been studied by Devroye and Fraiman \cite{Devroye}. The model presented in this paper is related to a case that has not been covered in \cite{Devroye} (see Section \ref{sec:devroye} for details).

\bigskip
\begin{teo}\label{connectedness}
Let $r\in \bbN\cup\{\infty\}$ be either a constant independent of $N$ or a function in $N$ such that $r=r(N)$ tends to a constant as $N\to \infty$. Let $i_*$ be defined by \eqref{maxi_*}, let $c(N)$ be a function in $N$ satisfying 
\begin{equation}\label{c(N)}
\quad-\ln k_{i_{*}} < c(N) < i_*(n-i_*)-\ln k_{i_{*}}
\end{equation} 
and let
\begin{equation}\label{c_*}
c_*:= \lim_{N\to \infty}c(N)\in[-\infty,\infty].
\end{equation} 
Consider the random graph $G(N,\K^r,p)$ with 
\begin{equation}\label{conn_threshold}
p:=\frac{\ln k_{i_{*}}+c(N)}{i_{*}(n-i_*)},
\end{equation}
where the condition \eqref{convki} holds.  
\begin{enumerate}
\item[(1)] If  $c_{*}=- \infty$, then 
$$\lim_{N\rightarrow\infty} \bbP[\,G(N,\K^r,p) \textit{ is connected }]=0.$$
\vspace{.1cm}
\item[(2)] If $c_{*}= c\in\bbR$ is a constant, then $$\lim_{N\rightarrow\infty} \bbP[\,G(N,\K^r,p) \textit{ is connected }] = e^{- e^{-c}-\gamma(i_*)},$$
where $\gamma(i_*)$ is defined as
\begin{equation}\label{gamma}
\gamma(i_*):= \lim_{N\to \infty} \sum_{i\neq i_*} k_i\, (k_{i_*}e^{c})^{-\frac{i(n-i)}{i_{*}(n-i_*)}}.
\end{equation}
\vspace{.1cm}
\item[(3)] If  $c_{*}=+ \infty$, then 
$$\lim_{N\rightarrow\infty} \bbP[\,G(N,\K^r,p) \textit{ is connected }]=1.$$
\end{enumerate}
\end{teo}

\vspace{.2cm}

In order to state our next results we need some more notation.  We  use the standard notation $\mu$-a.e. on some metric space $(\mathcal{S},\mu)$, which means that the set of elements for which a property does not hold is a set of $\mu$-measure zero. We shall suppress the measure $\mu$ in our notation if it is clear from context.  We use also the following standard asymptotic notation, in which the limits are taken as $N \rightarrow \infty$. For  functions $f =f(N)$ and $g=g(N)$, we write $f =O(g)$ if the limit of $f/g$ is bounded; $f = \Theta(g)$ if $f = O(g)$ and $g = O(f)$; $f = o(g)$ if $f/g \rightarrow 0$; $f\sim g$ if $f/g \rightarrow 1$.

Furthermore we denote the sizes of components of a graph $G$ by $L_1(G) \geq L_2(G) \geq \dots$, with $L_j(G) = 0$ if $G$ has fewer than $j$ components. The following result deals with the threshold for the existence of the giant component (composed of super-vertices) in $G(N,\K^r,p)$.

\bigskip 
\begin{prop}\label{corBollobas1}
Let $c\in \mathbb{R}^{+}$ be a positive constant and let $r\in \bbN\cup\{\infty\}$ be either a constant independent of $N$ or a function in $N$ such that $r=r(N)$ tends to a constant or $\infty$ as $N\rightarrow \infty$. 
Assume that  the limit $u:=\lim_{N\rightarrow\infty}(n/N)$ exists. Define 
$$s_2:=\sum_{i=1}^{r}\frac{i^2\, k_i}{n} \quad \text{and}\quad \bar{s_2}:=\lim_{N\to \infty} s_2,$$
provided the limit exists. Then the random graph  $G(N,\K^r,p)$ with $p=c/n$ satisfies the following properties.\\  

\begin{enumerate}
\item[(1)] If  $c\bar{s_2}\leq1$, then
$$\lim_{N\to \infty} \frac{L_1\left(G(N,\K^r,p)\right)}{N}=0 \quad\text{in probability},$$
while if $c\bar{s_2}>1$, then 
$$\lim_{N\to \infty}\bbP[\,L_1\left(G(N,\K^r,p)\right)=\Theta(N)\,]=1.$$

\item[(2)]  Let $\mu=(\mu_i)_{i\geq 1}$ be the probability measure on $\mathbb{Z}^{+}$ defined by \eqref{convki}. Then,
$$\lim_{N\to \infty}\,\frac{L_1\left(G(N,\K^r,p)\right)}{N}= \rho \quad\text{in probability},$$ 
where $\rho:=\sum_{i=1}^{r}h(i)\mu_i$ and the function $h$ is $\mu$-a.e. equal to the maximum solution of the non-linear equation $h=1-e^{-Th},$ where $T$ is an integral operator defined by $$Th(i):=\frac{ci}{u}\sum_{j=1}^{r}jh(j)\mu_j.$$ Furthermore, $\rho>0$ if and only if $c\bar{s_2}>1$.\\
\end{enumerate}
\end{prop}

Observe that $s_2$ is the average size of a super-vertex containing a random vertex. If the vertices in each super-vertex {\em were} all connected in $G(n,p)$, $s_2$ would be the {\em susceptibility} which is defined as the average component size of $G(n,p)$ (see e.g. \cite{JansonSpencer}). In Proposition \ref{corBollobas1} we assume that   the limit $\lim_{N\rightarrow\infty}(n/N)$ exists. This assumption guarantees that the sequence of kernels which will be given by (\ref{kernel}) is graphical in the sense of \cite{BJR} (see Sections \ref{relatedwork} and \ref{graphicalproof} for details). It is natural to ask what happens if the limit $\lim_{N\rightarrow\infty}(n/N)$ does not exist. More specifically, how does it affect the emergence and the size of the giant component in $G(N,\K^r,p)$?\\

In the next result we characterize the degree distribution in $G(N,\K^r,p)$. We define the degree of a super-vertex as the number of  super-vertices  connected to it in $G(N,\K^r,p)$.

\bigskip 
\begin{prop}\label{corBollobas}
Let $c\in \mathbb{R}^{+}$ be a positive constant and let $r\in \bbN\cup\{\infty\}$ be either a constant independent of $N$ or a function in $N$ such that $r=r(N)$ tends to a constant or $\infty$ as $N\rightarrow \infty$.
Assume that  the limit $u:=\lim_{N\rightarrow\infty}(n/N)$ exists. 
Consider the random graph $G(N,\K^r,p)$ with $p=c/n$. For each integer $k\geq 0$,  we let $Z_k$ denote the number  of super-vertices of degree $k$ in $G(N,\K^r,p)$. Then we have
$$\lim_{N\to \infty}\frac{Z_k}{N}= \bbP(\Xi = k) \quad\text{in probability},$$
where
\begin{equation}\label{mixedpoisson}
\bbP(\Xi = k):=\sum_{i\in\mathbb{Z^+}}\mu_i\, \bbP\left(Po(ci)=k\right)
\end{equation}
and $Po(ci)$ denotes a random variable with Poisson distribution with mean $ci$.
\end{prop}

\bigskip
The example below shows how and under which conditions the upper  tail of the degree distribution of a random super-vertex in $G(N,\K^{\infty},p)$  exhibits a power law behavior. In recent years, it has been conjectured that power laws characterize the behavior of the upper tails of the degree distribution in many real-worlds networks (see e.g. \cite{newman} for a review of the empirical evidence of this property). 

\bigskip 
\begin{exa}\label{powerlaw}
Consider the random graph  $G(N,\K^\infty,p)$ with $p=1/n$ and a configuration of super-vertices $\K^\infty$ such that
for some positive constants $C$ and $\alpha$,
$$
\sum_{i=k}^{\infty}\mu_{i} \sim \frac{C}{k^{\alpha}}
$$
for $k$ sufficient large. 
For each integer $k\geq 0$, we let $Z_{\geq k}$ denote the number of super-vertices of degree of at least $k$ in $G(N,\K^r,p)$. It follows as a corollary of Proposition \ref{corBollobas} that
\begin{equation}\label{example_powerlaw}
\frac{Z_{\geq k}}{N}\rightarrow \bbP(\Xi \geq k) \sim \frac{C}{k^{\alpha}},
\end{equation}

 where the first convergence is in probability when $N\rightarrow \infty$ (for $k$ fixed) and the latter approximation holds for $k$ large enough. To see it, note that \eqref{mixedpoisson} says that $\Xi$ follows a mixed Poisson distribution with mixing distribution. More precisely, $\Xi$ has a Poisson distribution with a random mean $Y$ such that $P(Y=i)=\mu_i$  for $i\ge 1$. 
We take $\epsilon>0$ arbitrarily small and obtain from \eqref{mixedpoisson}  that
$$\bbP(\Xi \geq k) \sim \bbP\left(\Xi \geq k | Y>(1-\epsilon)k\right) \frac{C}{k^{\alpha}} + \bbP\left(\Xi \geq k | Y<(1-\epsilon)k\right) \left(1 -\frac{C}{k^{\alpha}}\right).$$
Then \eqref{example_powerlaw}  follows from 
$$\bbP(\Xi \geq k | Y>(1-\epsilon)k) =1-o(1),$$
$$\bbP(\Xi \geq k | Y<(1-\epsilon)k) =o(k^{-\alpha}),$$
which may be obtained by Chernoff estimates (e.g. Corollary 13.1 in \cite{BJR}).  
\end{exa}

\bigskip 
This example shows that our model provides a mechanism starting from the critical  Erd\H os-R\'enyi random graph $G(n,1/n)$ and leading to scale-free networks. Therefore it may be seen as an instance of the good-get-richer mechanism proposed by Caldarelli et al. in \cite{caldarelli},  which has been introduced as an alternative to the well-known scale-free networks obtained by dynamical properties or preferential attachment. In \cite{caldarelli}, the authors start with a random graph with a large number $N$ of vertices. They associate  each vertex with the so-called fitness, which is a random number taken from a given probability distribution and measures the importance of the vertex. The probability of two vertices being connected  is given by a function of the fitnesses of the two vertices. If we assume that the fitness of a super-vertex in our model is its size that is distributed according to $(\mu_i)_{i\geq 1}$ defined by $\eqref{convki}$, then $G(N,\K^r,p)$ provides a {\em constructive} example of the good-get-richer mechanism.\\

\subsection{Related work}\label{relatedwork}

The random graph $G(N,\K^r,p)$ can be seen as a special case of inhomogeneous random graphs (IRG) studied by Bollob\'as, Janson and Riordan \cite{BJR}. To this end, we briefly recall some definitions and notations from \cite{BJR}. Consider a graph with vertex set $[N]:=\{1,2,\ldots,N\}$. A vertex space $\mathcal{V}$ is a triple $({\mathcal S},\mu,({\bf x}_N)_{N\geq 1})$ where ${\mathcal S}$ is a separable metric space, $\mu$ is a Borel probability measure on ${\mathcal S}$ and for each $N\in \mathbb N$, ${\bf x}_N$ is a random sequence $(x_1,x_2,\ldots,x_N)$ of points of ${\mathcal S}$ such that
\begin{equation}\label{condi}
\frac{|\{k:x_k \in A\}|}{N} \rightarrow \mu(A) \quad\text{in probability}\quad
\end{equation}
for every $\mu$-continuous set $A\subset {\mathcal S}$, where $|\mathcal{C}|$ denotes the cardinality of the set $\mathcal{C}$. A kernel $\kappa_N$ on  vertex space $\mathcal{V}$ is a symmetric non-negative Borel measurable function on ${\mathcal S}\times {\mathcal S}$.

Let $G^{\mathcal V}(N,\kappa_N)$ be the inhomogeneous random graph with vertex set $[N]$, in which two vertices $k$ and $l$ are connected by an edge  with probability
$$p_{kl}:=\min\left\{1,\frac{\kappa_N(x_k,x_l)}{N}\right\}.$$
This model is an extension of the one defined by S\"oderberg \cite{Sod} and its various properties are studied by \cite{BJR} under specific restrictions.  The main results regarding the existence and uniqueness of the giant component are proved by using an appropriated multi-type branching process and an integral operator  to which the component structure is related.

We shall show that our random graph model $G(N,\K^r,p)$  is a particular case of $G^{\mathcal V}(N,\kappa_N)$. To this end, we consider the set ${\mathcal S}=\mathbb{Z}^{+}$, the probability measure on ${\mathcal S}$ defined by $\mu(\{i\})=\mu_i$ given by (\ref{convki}), and the sequence ${\bf x}_N=(x_1,x_2,\ldots,x_N)$ of points of ${\mathcal S}$ such that $x_k$ represents the size of the $k$-$th$ super-vertex, for $k=1,2,\ldots,N$. Then the triple ${\mathcal V}:=({\mathcal S},\mu,({\bf x}_N)_{N\geq 1})$ is a vertex space. Observe that for all $i=1,2,\ldots,r$,
$$k_i=\sum_{k=1}^{N}I_{\{x_k = i\}},$$
where $I_A$ denotes the indicator random variable of the event $A$. Then by our construction the $x_k$'s are deterministic and therefore \eqref{convki} implies \eqref{condi}. Now we define the kernel $\kappa_N$ on the vertex space ${\mathcal V}$ by
\begin{equation}\label{kernel}
\kappa_N(x_k,x_l):=N\left(1-\left(1-p\right)^{x_kx_l}\right)
\end{equation}
and let the connection probabilities between two super-vertices of sizes $k$ and $l$ in  $G(N,\K^r,p)$  be defined as
$$p_{kl}:=\frac{\kappa_N(x_k,x_l)}{N}.$$ 
Then our model $G(N,\K,p)$  corresponds to $G^{\mathcal V}(N,\kappa_N)$.

%%%%%%%%%%%%%%%%%%%%%%%%%%%%%%%%%%%%%%%%%%%%%%%%%%%%%%%%%%%%%%%%%%%%%%%%%%%%%%
%%%%%%%%%%%%% PROOFS
%%%%%%%%%%%%%%%%%%%%%%%%%%%%%%%%%%%%%%%%%%%%%%%%%%%%%%%%%%%%%%%%%%%%%%%%%%%%%%

\bigskip
\section{Proof of Theorem \ref{connectedness}}\label{proofcon}
Throughout this section we let $r$ be  either a constant independent of $N$ or a function in $N$ such that $r=r(N)$ tends to a constant as $N\rightarrow \infty$. 
We may assume without loss of generality that $1\le r\le r^*$  for some constant $r^*\in \mathbb N$ (for $N$ sufficiently large).
We recollect some important conditions in Theorem \ref{connectedness}. We have
\begin{equation}\label{pconn2}
p:=\frac{\ln k_{i_{*}}+c(N)}{i_{*}(n-i_*)},
\end{equation} 
where $c(N)$ satisfies \eqref{c(N)} which ensures $p\in (0,1)$, and $i_{*}$ is defined by \eqref{maxi_*}, i.e. $i_{*}:=\Big\{ 1\le i \le r \ \Big|\  k_{i}^{\frac{1}{i (n-i)}} = \max_{1\le j\le r} k_{j}^{\frac{1}{j (n-j)}}\Big\}$.

Because the properties of being connected or of having no isolated super-vertices are monotone decreasing properties, we may assume that $c(N)=o(\ln N)$.

The connectedness of $G(N,\K^r,p)$ is closely related to the distribution of the number of isolated  super-vertices and ``small'' components, which will be the topics of the next three of sections and will be used in the proof of  Theorem \ref{connectedness} (1)--(3).

\bigskip
Let  $X$ denote the number of isolated super-vertices in $G(N,\K^r,p)$.
Observe that the random variable $X$ depends on $N$. For simplicity we suppress this dependence in our notation. 
First we shall derive asymptotic expressions of the first moment and second moment of $X$, which will be used in the proofs of Lemma \ref{zeroone-iso} (when $\lim_{N\to \infty}c(N) = \pm \infty$) and Lemma \ref{poisson} (when $\lim_{N\to \infty}c(N) \to c$ a constant).

\smallskip
Take any arbitrary  order of the super-vertices and let us write $X$ as a sum of indicator random variables
\begin{eqnarray}
X=\sum_{i=1}^r\sum_{k=1}^{k_i}I_k^i,\label{isolatedvertices}
\end{eqnarray}
where $I_k^i=1$ if  the $k$-$th$ super-vertex of size $i$ is isolated and 0 otherwise, for $k=1,2,\ldots,k_i$ and $i=1,2,\ldots,r$. 

Note that  $I_k^i=1$ if the $k$-th super-vertex of size $i$ is not connected with any other  super-vertex of size $i$ and  connected with any other super-vertex of size $j\neq i$ neither.  Since the super-vertices are connected independently of each other, we have
$$\bbE[I_k^i]=(1-p_{ii})^{k_i-1}\prod_{j\neq i}^r (1-p_{ij})^{k_j}=(1-p)^{i^2(k_i-1)}\prod_{j\neq i}^r (1-p)^{ijk_j}=(1-p)^{i(n-i)},$$
and hence
\begin{equation}\label{exp}
\bbE[X]=\sum_{i=1}^r\sum_{k=1}^{k_i}\bbE[I_k^i] = \sum_{i=1}^r k_i (1-p)^{i(n-i)}.
\end{equation}
Due to our choices of parameters, we have $N=\sum_{i=1}^{{r}}k_{i}\le n=\sum_{i=1}^{r}ik_{i}\le r N$, $c(N)=o(\ln N)$  and $p :=\frac{\ln k_{i_{*}}+c(N)}{i_{*}(n-i_*)}=o(1)$. Using these we obtain
\begin{eqnarray}\label{(1-p)power}
(1-p)^{i(n- i)}&=&\exp\left(-i(n- i) (p+O(p^2))\right)\nonumber\\
&=&\left(1+O\left(\frac{(\ln k_{i_*}+c(N))^2}{n}\right)\right)\, \exp\left(-\frac{i(n- i)}{i_{*}(n-i_*)} (\ln k_{i_*}+c(N))\right)\nonumber\\
%&=&\left(1+o(1)\right)\, \exp\left(-\frac{i(n- i)}{i_{*}(n-i_*)} (\ln k_{i_*}+c(N))\right)\nonumber\\
&=&\left(1+o(1)\right)\, k_{i_*}^{-\frac{i(n-i)}{i_{*}(n-i_*)}}\exp\left(-c(N) \frac{i(n-i)}{i_{*}(n-i_*)}\right).
\end{eqnarray}
To ease notation, we let
\begin{align}\label{funf(i)}
f(i):=f(i,c(N))
&= k_{i_*}^{\frac{i(n-i)}{i_{*}(n-i_*)}}\exp\left(c(N) \frac{i(n-i)}{i_{*}(n-i_*)}\right)
\end{align}
to obtain from \eqref{(1-p)power} that 
\begin{eqnarray}\label{(1-p)power2}
(1-p)^{i(n- i)}
&=&\left(1+o(1)\right)\,f(i)^{-1}\quad \text{for any} \quad 1\le i\le r.
\end{eqnarray}
By \eqref{exp}--\eqref{(1-p)power2}, the first moment of the number $X$ of isolated super-vertices in $G(N,\K^r,p)$ satisfies
\begin{align}\label{espe2}
\bbE[X]=\left(1+o(1)\right) \sum_{i=1}^r\frac{k_i}{f(i)}.
\end{align}

As for the second moment of $X$ we observe that
\begin{eqnarray}\label{secmom} 
\bbE[X^2]&=&\bbE\left[\left(\sum_{i=1}^r\sum_{k=1}^{k_i} I_k^i\right)^2\right]
=\bbE[X]-\sum_{i=1}^rk_i\bbE[I_1^iI_2^i]+\sum_{i,j=1}^{r}k_ik_j\bbE[I_1^iI_2^j].
\end{eqnarray}
Furthermore, we obtain by \eqref{(1-p)power2} that for any $1\le i,j\le r$
\begin{eqnarray}\label{exp2}
\bbE[I_1^iI_2^j]&=&\bbP(I_1^i=1\mid I_2^j=1)\bbP(I_2^j=1)\nonumber\\
&=&\left((1-p_{ii})^{k_i-1}(1-p_{ij})^{k_j-1}\prod_{l=1,l\neq i,j}^r(1-p_{il})^{k_l}\right)\bbP(I_2^j=1)\nonumber\\
&=&(1-p)^{i(n-i)+j(n-j)-ij}.
\end{eqnarray}
Note that  $(1-p)^{ij}=1+o(1)$ for any $1\le i,j\le r$. 
Thus  \eqref{secmom}--\eqref{exp2} and \eqref{(1-p)power2}--\eqref{espe2} imply
\begin{align} \label{second_moment2}
\bbE[X^2]&=\bbE[X]-\left(1+o(1)\right) \sum_{i=1}^r \frac{k_i}{f(i)^2} +\left(1+o(1)\right)  \sum_{i,j=1}^r \frac{k_i k_j}{f(i)f(j)}\nonumber\\
&=\bbE[X]-\left(1+o(1)\right) \sum_{i=1}^r \frac{k_i}{f(i)^2} +\left(1+o(1)\right) \left(\bbE[X]\right)^2.
\end{align}

\smallskip
Using \eqref{espe2} and \eqref{second_moment2}  we shall derive the existence of isolated super-vertices (in Lemma \ref{zeroone-iso}) and the exact distribution of the number of isolated super-vertices (in Lemma \ref{poisson}). 

\bigskip

\subsection{Existence of isolated super-vertices}
The following result deals with the existence of isolated super-vertices when $c_{*}=\pm  \infty$.
\smallskip
\begin{lem}\label{zeroone-iso} 
Let $r\in \bbN\cup\{\infty\}$ be either a constant independent of $N$ or a function in $N$ such that $r=r(N)$ tends to a constant as $N\to \infty$.  
Let  $X$ denote the number of isolated super-vertices in $G(N,\K^r,p)$ and let  
$c_*:= \lim_{N\to \infty}c(N)$ as in Theorem \ref{connectedness}. 
 
\begin{enumerate}
\item[(1)] 
If  $c_{*}=+ \infty$, then 
$$
\lim_{N\rightarrow\infty} \bbP[X\ge 1]=0.$$
\item[(2)] If  $c_{*}=- \infty$, then 
$$\lim_{N\rightarrow\infty} \bbP[X\ge 1]=1.$$
\end{enumerate}
%$$
%\lim_{N\rightarrow\infty} \bbP[X\ge 1]=
%\begin{cases} 
%0, \quad \text{if}\quad c_{*}=+ \infty,\\
%\end{cases}
%$$
\end{lem}

\smallskip
\begin{proof}
Without loss of generality we may assume that $1\le r\le r^*$  for some constant $r^*\in \mathbb N$ (for $N$ sufficiently large).

\smallskip
To prove (1), we assume  without loss of generality that $c(N)>0$ (for some large $N$). 
Note that $\frac{i(n-i)}{i_{*}(n-i_*)}\ge \frac{(n-i)}{i_{*}(n-i_*)}=\frac{1+o(1)}{i_{*}}$ and thus 
\begin{eqnarray}\label{reciprocal-fi}
f(i)^{-1}&=&k_{i_*}^{-\frac{i(n-i)}{i_{*}(n-i_*)}}\exp\left(-c(N) \frac{i(n-i)}{i_{*}(n-i_*)}\right)\nonumber\\
&\le &k_{i_*}^{-\frac{i(n-i)}{i_{*}(n-i_*)}} \exp\left(- \left(1+o(1)\right)\frac{1}{i_{*}}c(N)\right).
\end{eqnarray}
Furthermore, by definition of $i_*:=\Big\{ 1\le i \le r \ \Big|\  k_{i}^{\frac{1}{i (n-i)}} = \max_{1\le j\le r} k_{j}^{\frac{1}{j (n-j)}}\Big\}$, we have that for each $1\le i\le r$
\begin{eqnarray}\label{ki}
k_i \leq k_{i_*}^{\frac{i(n-i)}{i_*(n-i_*)}},
\end{eqnarray}
which together with \eqref{reciprocal-fi} implies that for each  $1\le i\le r$
\begin{align*}%\label{kivsf(i)} 
\frac{k_i}{f(i)} 
&\leq \exp\left(- \left(1+o(1)\right)\frac{1}{i_{*}}c(N)\right).
\end{align*}
From this in \eqref{espe2}, we obtain
\begin{eqnarray*}
\bbE[X] &=& \left(1+o(1)\right) \sum_{i=1}^r\frac{k_i}{f(i)}\\  
 &\leq&    \left(1+o(1)\right)r\exp\left(- \left(1+o(1)\right)\frac{1}{i_{*}}c(N)\right)\\
&\leq&\left(1+o(1)\right) r^*\exp\left(- \left(1+o(1)\right)\frac{1}{r_{*}}c(N)\right). 
\end{eqnarray*}
(The last inequality follows, because $r \le r^*$ and $ i_{*} \le r^*$.)
Thus by Markov's inequality, we have $\bbP(X\ge 1)\le \bbE[X]  \to 0 $ as $N\to \infty$, when $\lim_{N\to \infty}c(N)=+ \infty$. 

\bigskip
To prove (2), we  observe that
\begin{align*}%\label{kivsf(i)} 
\frac{k_{i_*}}{f(i_*)} 
&= k_{i_*} k_{i_*}^{-\frac{i_*(n-i_*)}{i_{*}(n-i_*)}}\exp\left(-c(N) \frac{i_*(n-i_*)}{i_{*}(n-i_*)}\right) 
=\exp\left(-c(N) \right)\\
\end{align*}
and since $\sum_{i=1}^r\frac{k_i}{f(i)}\geq \frac{k_{i_*}}{f(i_*)}$ and $c(N)\to -\infty$ as $N\to \infty $, we have
 \begin{align}
\bbE[X]&=\left(1+o(1)\right) \sum_{i=1}^r\frac{k_i}{f(i)}\nonumber\\
&\geq  \left(1+o(1)\right) \frac{k_{i_*}}{f(i_*)} =  \left(1+o(1)\right) \exp\left(-c(N) \right) \to+ \infty \quad\text{as}\quad N\to  \infty.\label{espe5}
\end{align}
We shall prove that $\bbP(X=0)\to 0$ as $N\to \infty$, 
by applying Chebyshev's inequality $$\bbP(X=0)\le \frac{\bbE[X^2]}{(\bbE[X])^2}-1.$$  
It suffices to show that  $\lim_{N\to \infty}\frac{\bbE[X^2]}{(\bbE[X])^2} = 1$ as $N\to \infty$.  
To this end, we observe that \eqref{second_moment2} implies
\begin{align}\label{second_moment5}
\frac{\bbE[X^2]}{(\bbE[X])^2}
&=\bbE[X]^{-1}-\left(1+o(1)\right) \sum_{i=1}^r \frac{k_i}{f(i)^2} \left(\bbE[X]\right)^{-2}  +1+o(1).
\end{align}
Using the property that $f(i)>1$ for each $1\le i\le r$ (because $c(N)>-\ln k_{i_*}$), we have
\begin{align*}%\label{comparison2}
 \sum_{i=1}^r\frac{k_i}{f(i)^{2}}\le \sum_{i=1}^r\frac{k_i}{f(i)}=\left(1+o(1) \right) \bbE[X].
\end{align*}  
Since $\bbE[X]\to+ \infty$  as $N\to \infty $ by \eqref{espe5}, we obtain 
\begin{align}\label{comparison3}
 \left(\sum_{i=1}^r\frac{k_i}{f(i)^{2}}\right)   \left(\bbE[X]\right)^{-2}  
 &< \left(1+o(1) \right) (\bbE[X] )^{-1} \to 0\quad\text{as}\quad N\to \infty.
\end{align}
Putting \eqref{second_moment5}  and \eqref{comparison3} together, we have
$\frac{\bbE[X^2]}{(\bbE[X])^2}  \to 1$ as $N\to \infty$ as desired.
\end{proof}

\bigskip
\subsection{Distribution of the number of isolated super-vertices}
In this section we deal with the exact asymptotic distribution of the number of isolated super-vertices when $c_* =c\in \mathbb{R}$ is a constant, by applying  Stein's method, in particular Theorem 6.24 \cite{JLR}. To this end we need a few more definitions and notations.  The random variables $\{I_k\}_k$ are said to be positively related if they satisfy the following two conditions: 
\begin{enumerate}
\item[(1)]  For each $k$ there exists a family of random variables $J_{\ell}^k$, $\ell\neq k$, such that the joint distribution of $\{J_{\ell}^k\}_{\ell}$ is the same as the conditional distribution of $\{I_{\ell}\}_{\ell}$ given $I_k=1$; and 
\item[(2)]   $J_{\ell}^k\geq I_{\ell}$ for every $\ell\neq k$.
\end{enumerate} 
\smallskip
Formally, 
$\mathcal{L}(\{J_{\ell}^k\}_{\ell})=\mathcal{L}(\{I_{\ell}\}_{\ell}\mid I_k=1),$ and $J_{\ell}^k\geq I_{\ell}$ for every ${\ell}\neq k$, 
where $\mathcal{L}(\{Y_{\ell}\}_{\ell})$ denotes the joint distribution of the random variables $\{Y_{\ell}\}_{\ell}$.

\bigskip
\begin{lem}[Theorem 6.24 in \cite{JLR}]\label{totaldev}
Given random variables $\{I_k\}_k$, let $X:=\sum_{k} I_k$.  
If $\{I_k\}_k$ are positively related, then
\begin{equation}\label{vardist}
d_{TV}(X,Po(\lambda))\leq\frac{\bbV(X)}{\bbE[X]}-1+2\max_{1\leq k\leq N}\{\bbE(I_k)\},
\end{equation}
where $\lambda:=\bbE[X]$ and $d_{TV}(\cdot,\cdot)$ is the  total variation distance.
\end{lem}

\bigskip
\begin{lem}\label{poisson}
Let $r\in \bbN\cup\{\infty\}$ be either a constant independent of $N$ or a function in $N$ such that $r=r(N)$ tends to a constant as $N\to \infty$. 
Let  $X$ denote the number of isolated super-vertices in $G(N,\K^r,p)$.
If $c_{*}=c\in\bbR$ is a constant, then $X$ has asymptotically Poisson distribution with mean  $e^{-c}+\gamma(i_*)$, 
where $\gamma(i_*):= \lim_{N\to \infty} \sum_{i\neq i_*} k_i\, (k_{i_*}e^{c})^{-\frac{i(n-i)}{i_{*}(n-i_*)}}$ is defined as in \eqref{gamma}.
\end{lem}

\smallskip
\begin{proof}
Without loss of generality we may assume that $1\le r\le r^*$  for some constant $r^*\in \mathbb N$ (for $N$ sufficiently large) and that $c(N) = c\in \bbR$ is a constant.

By the definition of $i^*$, we have that %$\frac{k_{i_*}}{f(i_*)}= e^{-c}$ and 
$k_i \leq k_{i_*}^{\frac{i(n-i)}{i_*(n-i_*)}}$ for each $1\le i\le r$ and $1\le i_{*}\le r$. 
Note further that 
$ (1+o(1))\frac{1}{i_{*}}\le \frac{i(n-i)}{i_{*}(n-i_*)}\le(1+o(1))\frac{r^*}{i_{*}}$ for each $1\le i\le r$. 
Thus, letting  $M(0):=0$, $M(c):=\left(-(1+o(1))\,  \frac{c}{ r^*}\right)$ if $c>0$ 
and $M(c):=\left(-(1+o(1))\, c\, r^*\right)$ if $c<0$ (note that $M(c)$ is independent of $N$), we have that for each   $1\le i\le r$
\begin{align*}%\label{M(c)}
c\, \frac{i(n-i)}{i_{*}(n-i_*)} \ge -M(c)
\end{align*}
and hence
\begin{align*}
\frac{k_i}{f(i)}= k_i k_{i_*}^{-\frac{i(n-i)}{i_{*}(n-i_*)}}\exp\left(-c \frac{i(n-i)}{i_{*}(n-i_*)}\right)
\ \le \ \exp(M(c)).
\end{align*}
Therefore, we have
\begin{align}\label{largei}
\gamma(i_*):=\lim_{N\to \infty} \sum_{i\neq i_*} k_i\, (k_{i_*}e^{c})^{-\frac{i(n-i)}{i_{*}(n-i_*)}}=\lim_{N\to \infty} \sum_{i\neq i_*} \frac{k_i}{f(i)} \ \le \  r^* \exp(M(c))< \infty.
\end{align} 
Because $\frac{k_{i_*}}{f(i_*)}= e^{-c}$, \eqref{espe2} and \eqref{largei}  imply that
\begin{align}\label{limitEx}
\bbE[X]
 %&=\left(1+o(1)\right) \sum_{i=1}^r\frac{k_i}{f(i)}\nonumber\\  
&= \left(1+o(1)\right)\left( \frac{k_{i_*}}{f(i_*)} +    \sum_{i\neq i_*}\frac{k_i}{f(i)} \right)\to e^{-c} +\gamma(i_*) \quad\text{as}\quad N\to \infty.
\end{align}

\smallskip
In order to apply \eqref{vardist}, we take an arbitrary  order of the super-vertices and rewrite the number $X$ of isolated super-vertices in $G(N,\K^r,p)$ as
\begin{equation}\label{newindicators}
X:=\sum_{k=1}^N I_k,
\end{equation}
where for each $k=1,2,\ldots,N$, $I_k=1$ if the $k$-th super-vertex is isolated and $I_k=0$ otherwise.  
Let $G_k(N,\K^r,p)$ be the random graph $G(N,\K^r,p)$ with all edges from the $k$-th super-vertex  removed,  and let $J_{\ell}^k=1$ if the $\ell$-th super-vertex is isolated in $G_k(N,\K^r,p)$, and $J_{\ell}^k=0$ otherwise. Observe that
$$J_{\ell}^k=(I_{\ell}\mid I_k=1),$$ and thus
$$\mathcal{L}(\{J_{\ell}^k\}_{\ell})=\mathcal{L}(\{I_{\ell}\}_{\ell}\mid I_k=1).$$ Moreover, if the $\ell$-th super-vertex in $G(N,\K^r,p)$ is isolated (i.e. $I_{\ell}=1$), then $J_{\ell}^k=1$ as the $\ell$-th super-vertex is still isolated even if all the edges from the $k$-th super-vertex are removed. Otherwise, i.e. $I_{\ell}=0$, then the $\ell$-th super-vertex is either connected with the $k$-th super-vertex or with any other super-vertex. So, if all edges from the $k$-th super-vertex are removed, the resulting $\ell$-th super-vertex could  become either isolated or not, this means that  $J_{\ell}^k=\{0,1\}$. Thus $J_{\ell}^k\geq I_{\ell}$ for every ${\ell}\neq k$, and thus the random variables $\{I_k\}_k$ are positively related. So we can apply Lemma \ref{totaldev}. 

Next we shall show that $d_{TV}(X,Po(\lambda))\to 0$  as $N\to \infty$, by proving $\frac{\bbV(X)}{\bbE[X]} \to 1$ and $\max_{1\leq k\leq N}\{\bbE(I_k)\} \to 0$ as $N\to \infty$.  
By \eqref{espe2} and \eqref{second_moment2} as well as \eqref{limitEx}, we get 
\begin{align}\label{second_moment8}
\frac{\bbV(X)}{\bbE[X]} &=  \frac{\bbE[X^2] - (\bbE[X])^2}{\bbE[X]}\nonumber\\
&=1 -\left(1+o(1)\right) \sum_{i=1}^r \frac{k_i}{f(i)^2} (\bbE[X]))^{-1}  + o(1)\bbE[X]\nonumber\\
&=1 + o(1) -\left(1+o(1)\right) \sum_{i=1}^r \frac{k_i}{f(i)^2} (\bbE[X]))^{-1}.  
\end{align}
Let $f_*:=\min\{f(i): i \in \{1,2,\ldots, r\}\}$.
By \eqref{espe2} we have 
\begin{align}\label{second_moment9}
\sum_{i=1}^r \frac{k_i}{f(i)^2}\le \sum_{i=1}^r \frac{k_i}{f(i) f_*} \le   (1+o(1)) \frac{\bbE[X]}{f_*}.
\end{align}  

By \eqref{second_moment8}, in order to prove 
\begin{align}\label{second_moment4}
\frac{\bbV(X)}{\bbE[X]} \to 1  \quad \text{as}\quad  N\to \infty, 
\end{align} 
it suffices to show that 
\begin{align}\label{todominf(i)}
f_* \to \infty \quad \text{as}\quad N\to \infty,
\end{align} 
because this together with   \eqref{second_moment9}  implies that
$$\sum_{i=1}^r \frac{k_i}{f(i)^2}  (\bbE[X]))^{-1} \le (1+o(1)) \frac{1}{f_*} 
\to   0 \quad \text{as}\quad N\to \infty.$$

In order to prove \eqref{todominf(i)}, we observe from the definition of $M(c)$ that for each $1\le i\le r$
\begin{align*}
c\, \frac{i(n-i)}{i_{*}(n-i_*)} \ge -M(c),
\end{align*}
and therefore we get
\begin{align}\label{lowerboundf(i)}
f(i)
&= k_{i_*}^{\frac{i(n-i)}{i_{*}(n-i_*)}}\exp\left(c\, \frac{i(n-i)}{i_{*}(n-i_*)}\right)\nonumber\\
&\ge  k_{i_*}^{\frac{i(n-i)}{i_{*}(n-i_*)}}\exp\left(-M(c)\right)\nonumber\\
&\ge   k_{i_*}^{ (1+o(1))\frac{1}{i_{*}}}\exp\left(-M(c)\right).
\end{align}
Now we define
\begin{align}\label{ell*}
\ell_*:=\min\{\ell \in\{1,2,\ldots,r\} : \mu_{\ell}>0\}.
 \end{align} 
By the definition of $i_*:=\Big\{ 1\le i \le r \ \Big|\  k_{i}^{\frac{1}{i (n-i)}} = \max_{1\le j\le r} k_{j}^{\frac{1}{j (n-j)}}\Big\}$, we have  $i_*\le \ell_*$ and thus  
$$k_{i_*} \geq k_{\ell_*}^{\frac{i_*(n-i_*)}{\ell_*(n-\ell_*)}}
= k_{\ell_*}^{ (1+o(1))\frac{i_*}{\ell_*}}.$$ 
Putting this in \eqref{lowerboundf(i)} we have that for each $1\le i\le r$
\begin{align*}
f(i) 
&\ge   k_{i_*}^{ (1+o(1))\frac{1}{i_{*}}}\exp\left(-M(c)\right)
\ge k_{\ell_*}^{ (1+o(1))\frac{1}{\ell_*}}\exp\left(-M(c)\right).
\end{align*}
Because $k_{\ell_*}=(1+o(1))\mu_{\ell_*}N$ with $\mu_{\ell_*}>0$, we get
\begin{align*}
f_* &:=\min_{1\le i\le r} f(i)\\
& \ge  \left( (1+o(1))\mu_{\ell_*}N\right)^{(1+o(1))\frac{1}{\ell_*}} \exp\left(-M(c)\right)
\to \infty \quad \text{as}\quad N\to \infty.
\end{align*}

\smallskip
Finally we shall show that $\max_{1\leq k\leq N}\{\bbE(I_k)\} =o(1)$. To this end, we consider the sequence of indicator random variables $(I_k)_k$ used in \eqref{newindicators}. For each $k=1,2,\ldots,N$  we obtain
$$\bbE(I_k)=\sum_{i=1}^{r}(1-p)^{i(n-i)}\frac{k_i}{N} = \frac{\bbE[X]}{N} \to 0 \quad\text{as}\quad N\to \infty,$$ 
where the last equality and the limit behavior follow  from \eqref{exp}  and \eqref{limitEx} respectively. 
Therefore we have
\begin{align}\label{maxIk}
 \max_{1\leq k\leq N}\{\bbE(I_k)\}=o(1).
\end{align}
Thus, we can conclude from \eqref{vardist}, \eqref{second_moment4} and \eqref{maxIk}   that  $X$ has asymptotically Poisson distribution with mean  $ e^{-c}+\gamma(i_*)$.
\end{proof}

\bigskip
\subsection{Small components}
Before we proceed to the proof of Theorem \ref{connectedness}, we finally consider components of size $m$ for each $2\leq m\leq N/2$. 

\smallskip
\begin{lem}\label{smallcomponents}
Let $r\in \bbN\cup\{\infty\}$ be either a constant independent of $N$ or a function in $N$ such that $r=r(N)$ tends to a constant as $N\to \infty$. 
If $c_{*}= + \infty$ or $c_* =c\in \mathbb{R}$ is a constant, then  $whp$ $G(N,\K^r,p)$ does not have any component of size $m$ for any $2\leq m\leq N/2$.
\end{lem}

\smallskip
\begin{proof}
For each $2\leq m \leq N/2$, we let $\mathcal{S}^m$ be the set of all subsets of $m$ super-vertices, and for each $S\in \mathcal{S}^m$, let $m_i(S)$ be the number  of super-vertices of size $i$ in $S$.  Note that  $m=m_1(S)+\dots+m_r(S)$. Observe that if the super-vertices in $S$ form a component in $G(N,\K^r,p)$, then the following two events should hold:

\begin{enumerate}
\item[(i)] $A^S_1:=\{\text{the super-vertices in }S \text{ are connected}\}$ and
\item[(ii)] $A^S_2:=\{\text{no super-vertex in } S \text{ is connected with a super-vertex in } S^c\}$,
\end{enumerate}
where $S^c$ denotes the complement of $S$. 

Note that the events  $A^S_1$ and $A^S_2$ are independent, because every pair of vertices are connected independently of each other with probability $p$. Thus, we have
\begin{eqnarray}\label{exi} \displaystyle
\bbP(\exists \text{ a component of size }m) \leq&\sum_{S\in \mathcal{S}^m}\bbP( A^S_1)\bbP( A^S_2).
\end{eqnarray}

Since a component of size $m$ contains a tree of size $m$ and the number of Cayley trees on $m$ vertices is $m^{m-2}$,  moreover $p_{ij}=1-\left(1-p\right)^{ij}\leq ijp\leq r^2p$,  we have
\begin{eqnarray}\label{A1}
\bbP(A^S_1)\leq \bbP(S\text{ contains a tree })\leq m^{m-2}(r^2p)^{m-1}.
\end{eqnarray}
On the other hand, for each $S\in \mathcal{S}^m$ letting $M_{m}(S):=\sum_{i=1}^r im_i(S)$,  we have
$$
 \bbP(A^S_2)=(1-p)^{M_{m}(S)(n-M_{m}(S))}.
 $$
Observe that $M_{m}(S) \geq m$. Since the function $f(x)=x(n-x)$ is increasing in $x\in (2,N/2)$, we have
\begin{equation}\label{A2} 
 \bbP(A^S_2)\leq(1-p)^{m(n-m)}.
 \end{equation}

By \eqref{exi}--\eqref{A2}, we obtain
\begin{align}
\label{exista}
\bbP(\exists \text{ a component of size }m) 
&\le& \binom{N}{m}  m^{m-2}(r^2p)^{m-1}(1-p)^{m(n-m)}\nonumber\\
&\le& (eN)^m\, m^{-5/2}\,(r^2 p)^{m-1}\, e^{-p(m-1)(n-m)},
\end{align}
where the last inequality is because $\binom{N}{m}<(eN)^m/m^{m+1/2}$ and $(1-p)< e^{-p}$.
Summing up, we have
\begin{align}
\label{existsum}
\sum_{m=2}^{N/2}\bbP(\exists \text{ a component of size }m) 
\le \sum_{m=2}^{N/2} (eN)^m\, m^{-5/2}\,(r^2 p)^{m-1}\, e^{-p(m-1)(n-m)}.
\end{align}

 We shall show that the left hand side of \eqref{existsum} tends to 0 as $N\to \infty$, making a case distinction depending on whether $ 2\ell_* +2\le m\le N/2$ or $1\le m\leq 2\ell_* +1$, where $\ell_*:=\min\{\ell \in\{1,2,\ldots,r\} : \mu_{\ell}>0\}$ as defined in \eqref{ell*}.\\  

\noindent
{\bf Case 1.} Let $2\ell_* +2\le m\le N/2$. 
We first rewrite \eqref{exista}  as 
\begin{equation}\label{eq1}
eN \, m^{-5/2}\, (er^2N)^{m-1}\,  \exp\left\{ (m-1) \left[\ln p - p\,(n-m)\right]\right\}.
\end{equation}
Note that because $p:=\frac{\ln k_{i_{*}}+c(N)}{i_{*}(n-i_*)}$,  we have
\begin{equation*}\label{eq2}
\ln p - p\,(n-m)=\ln\left[\ln k_{i_*} + c(N)\right] - \ln\left[ i_*(n-i_*)\right]- \left(\frac{\ln k_{i_*} + c(N)}{i_*}\right)\,\frac{(n-m)}{(n-i_*)},
\end{equation*}
so considering that  $m\leq N/2 < n/2$ (and thus $\frac{n-m}{n-i_*}\ge \frac{1}{2}$), we obtain
\begin{equation*}\label{eq3}
\ln p - p\,(n-m)\le \ln \left[\ln k_{i_*} + c(N)\right] - \frac{1}{2\, i_*} (\ln k_{i_*} + c(N)) -  \ln\left[ i_*(n-i_*)\right].
\end{equation*}
Observe now that for any $K>0$ there exists $\tilde x>0$ such that $\ln x<x/K$ for any $x>\tilde x$.  Using this with $x=\ln k_{i_*}+c(N)$, we have $\ln \left[\ln k_{i_*} + c(N)\right] < \frac{1}{K} (\ln k_{i_*} + c(N))$ and so we have
\begin{equation}\label{eq5}
\ln p - p\,(n-m)\le \left(\frac{1}{K}-\frac{1}{2i_*}\right)\left(\ln k_{i_*} + c(N)\right)-  \ln\left[ i_*(n-i_*)\right],
\end{equation}   
provided $N$ is sufficiently large. Using \eqref{eq5},  
we can bound (\ref{eq1}) from above by
\begin{align*}
&& e N  \, m^{-5/2}\, (er^2N)^{m-1}\,  \exp\left\{ (m-1) \left[ \left(\frac{1}{K}-\frac{1}{2i_*}\right)\left(\ln k_{i_*} + c(N)\right)-  \ln\left[ i_*(n-i_*)\right]\right]\right\}\\
&&=
 eN  \, m^{-5/2}\, \left(\frac{er^2N}{i_*(n-i_*)}\right)^{m-1} \, \exp\left\{ (m-1) \left[ \left(\frac{1}{K}-\frac{1}{2i_*}\right) \left(\ln k_{i_*}+ c(N)\right)\right]    \right\}\\
&&\le eN  \, m^{-5/2}\, \left(er^2\right)^{m-1} \, \exp\left\{ (m-1) \left[ \left(\frac{1}{K}-\frac{1}{2i_*}\right) \left(\ln k_{i_*}+ c(N)\right)\right]    \right\}\\
&&\le eN  \, m^{-5/2}\,   \exp\left\{ (m-1) \left[ \left(\frac{1}{K}-\frac{1}{2i_*}\right) \left(\ln k_{i_*}+ c(N)\right)+\ln(er^2)\right]    \right\}.
 \end{align*}
The last inequality follows because $N/i_*(n-i_*) \le 1$ for $N$ sufficiently large.  
Therefore, we have  that for 
\begin{align}\label{case11}
&& \bbP(\exists \text{ a component of size }m)\nonumber\\
&\leq &  eN  \, m^{-5/2}\, \exp\left\{ (m-1) \left[   \left(\frac{1}{K}-\frac{1}{2 i_*}\right) \left(\ln k_{i_*}+ c(N)\right)+\ln(er^2)\right]\right\}.
\end{align}
Now consider a constant $a>1$. Taking $K$ large enough such that $K > 2i_*(2\ell_* +1)a$  we have
\begin{equation*}
\begin{array}{llll}
&&\hspace{-1.3cm} (m-1)\left(\frac{1}{K}-\frac{1}{2i_*}\right)\\[.4cm]
&\hspace{-0.3cm} \leq &\hspace{-0.3cm} (2\ell_*+1)\left(\frac{1}{K}-\frac{1}{2i_*}\right)& \text{ (because }m\geq 2\ell_* +2)\\[.4cm]
&\hspace{-0.3cm}\leq &\hspace{-0.3cm}  \left(\frac{1-(2\ell_*+1)a }{2i_*a} \right) & \text{ (because } \frac{1}{K} < \frac{1}{2i_*(2\ell_* +1)a}). 
\end{array}
\end{equation*}
Thus we have
\begin{equation}\label{case12}
\begin{array}{llll}
&&\hspace{-1.3cm} (m-1)\left(\frac{1}{K}-\frac{1}{2i_*}\right)\ln k_{i_*}\\[.4cm]
&\hspace{-0.3cm}\leq &\hspace{-0.3cm}  \left(\frac{1-(2\ell_*+1)a }{2i_*a} \right) \ln k_{i_*}  =\left(\frac{(2\ell_*+1)a -1}{2i_*a} \right) \ln k_{i_*}^{-1}&  \\[.4cm]
&\hspace{-0.3cm}\leq &\hspace{-0.3cm} -(1+o(1))\left(1+\frac{a-1}{2a\ell_*}\right) \ln k_{\ell_*} & \text{ (because }k_{i_*}\geq k_{\ell_*}^{\frac{i_*(n-i_*)}{\ell_* (n-\ell_*)}}).
\end{array}
\end{equation}
Furthermore, because $m\geq 2\ell_* +2$, we have
\begin{equation}\label{case122}
\begin{array}{llll}
&&\hspace{-1.3cm} (m-1)\left[\left(\frac{1}{K}-\frac{1}{2 i_*}\right)c(N) + \ln\left(er^2\right)\right] \\[.4cm]
&\hspace{-0.3cm}\leq &\hspace{-0.3cm} -(2\ell_* + 1)\left[\left(\frac{1}{2 i_*}-\frac{1}{K}\right)c(N) - \ln\left(er^2\right)\right],
\end{array}
\end{equation}
and   because  $c(N)\to +\infty$ or $c\in \mathbb R$, we have
$$\exp\left\{1-(2\ell_* + 1)\left[\left(\frac{1}{2 i_*}-\frac{1}{K}\right)c(N) - \ln\left(er^2\right)\right]\right\}  \le C_1,$$
for $K$ and $N$ large enough.
Using \eqref{case11}--\eqref{case122}, we have  that for $2\ell_* +2\le m\le N/2$
\begin{align}\label{case13}
\bbP(\exists \text{ a component of size }m)
\, \leq \, 
C_1  \underbrace{\left(N k_{\ell_*}^{-(1+o(1))\left(1+\frac{a-1}{2a\ell_*}\right)}\right)}_{=:\varphi_1(N)} m^{-5/2}
\end{align}
where $\varphi_1(N) \rightarrow 0$ as $N\rightarrow \infty$, because $k_{\ell_*}=(1+o(1))\mu_{\ell_*}N$, $\mu_{\ell_*}>0$ and $a>1$.

Summing up, we have
\begin{equation}\label{caselargem}
\sum_{m=2\ell_{*}+2}^{N/2}\bbP(\exists \text{ a component of size $m$})
 \, \leq  \,  C_1 \varphi_1(N)\sum_{m=2\ell_{*}+2}^{N/2} m^{-5/2},
\end{equation}
which tends to 0 as $N\to \infty$.

\smallskip
\noindent
{\bf Case 2.} 
Let $2\le m\leq 2\ell_* +1$. In this case, it suffices to consider components containing only super-vertices of size at least $\ell_*:=\min\{\ell \in\{1,2,\ldots,r\} : \mu_{\ell}>0\}$.

To see this let $E$ be the event that there exists a component of size $m$ in which all its super-vertices have size at least $\ell_*$, and let $Y$ be the number of super-vertices of size strictly less than $\ell_*$ in a random sample of $m$ super-vertices. Hence, $Y$ is a random variable that follows the hypergeometric distribution with parameters $(N,v_{\ell_*},m)$, where $v_{\ell_*}$ denotes the number of super-vertices of size  strictly less than $\ell_*$, i.e., $v_{\ell_*}:=\sum_{i=1}^{\ell_*-1}k_i$.  Observe that, as $v_{\ell_*} =o(N)$ and $m\leq 2\ell_* +1$, so $m$ is independent of $N$, 
\begin{equation*}
\bbP(Y=0)=\frac{\binom{N-v_{\ell_*}}{m}}{\binom{N}{m}}=\frac{(N-v_{\ell_*})}{N} \cdots \frac{(N-v_{\ell_*}-m+1)}{(N-m+1)}=1-o(1),
\end{equation*}
and therefore
\begin{equation*}
\bbP(E^c)\leq \bbP(Y>0)=o(1).
\end{equation*}
 
More precisely, we obtain
\begin{equation}\label{existmsmall}
\bbP(\exists \text{ a component of size }m)= \bbP(E)+ o(1).
 \end{equation}
Thus, applying  arguments analogous to \eqref{exi}-\eqref{exista} to $\bbP(E)$, and noting that in such case the component $S$ has only super-vertices of size at least $\ell_*$, we obtain that  $M_m(S)\geq m\ell_*$  and therefore, instead of \eqref{exista}, we obtain
\begin{eqnarray}\label{caso2otravez}
\bbP(E) 
&\le& (eN)^m\, m^{-5/2}\,(r^2 p)^{m-1}\, e^{-p(m\ell_*-1)(n-m\ell_*)}\nonumber\\
&=&O(N)\, (eNp)^{m-1}\, e^{-p m\ell_*(n-m\ell_*)+m\ln r^2}\quad (\text{ because } m\ge 2\text{ and } N\le n).
\end{eqnarray}
Because $2\leq m\leq2\ell_*+1$, we have $m\ell_*(n-m\ell_*)\ge 2\ell_*[n-\ell_*(2\ell_*+1)]$. 
Using $p:=\frac{\ln k_{i_{*}}+c(N)}{i_{*}(n-i_*)}$ and $k_{i_*}\geq k_{\ell_*}^{\frac{i_*(n-i_*)}{\ell_*(n-\ell_*)}}$, we therefore obtain
\begin{eqnarray*}
p\ell_* m\,(n-\ell_* m) &\geq& (\ln k_{i_*} + c(N))\frac{2\ell_*[n-\ell_*(2\ell_*+1)]}{i_*(n-i_*)}\\
&\geq&  (\ln k_{\ell_*})\frac{2[n-\ell_*(2\ell_*+1)]}{(n-\ell_*)}+c(N)\frac{2\ell_*[n-\ell_*(2\ell_*+1)]}{i_*(n-i_*)}\\
&=&(2+o(1)) (\ln k_{\ell_*})+(2+o(1))\frac{\ell_*}{i_*}c(N),
\end{eqnarray*}
and  
\begin{eqnarray*}
(eNp)^{m-1}\le \left(\frac{eN[\ln k_{i_{*}}+c(N)]}{i_{*}(n-i_*)}\right)^{2\ell_*}
=O(1)[\ln k_{i_{*}}+c(N)]^{2\ell_*} \quad (\text{ because } N\le n).
\end{eqnarray*}
Putting these  in (\ref{caso2otravez}), because $m\le N/2$ we obtain
\begin{eqnarray*}\label{lastequa}
\bbP(E) 
&\le&O(N)\, [\ln k_{i_{*}}+c(N)]^{2\ell_*}\,  \left[k_{\ell_*}\right]^{-2+o(1)}\, e^{-(2+o(1))\frac{\ell_*}{i_*}c(N)+ (2\ell_*+1)\ln r^2} \\
&=:&\varphi_2(N) \to 0 \quad \text{ as } N \to \infty, 
\end{eqnarray*}
because the condition $c(N)\to +\infty$ or $c\in \mathbb R$ implies $e^{-(2+o(1))\frac{\ell_*}{i_*}c(N)+ (2\ell_*+1)\ln r^2}=O(1)$ and  the conditions $k_{\ell_*} =(1+o(1))\mu_{\ell_*}N$ with $\mu_{\ell_*}>0$ and $k_{i_*}\leq O(N)$ imply  $O(N)\, (\ln k_{i_*}+c(N))^{2\ell_*} \,  \left[k_{\ell_*}\right]^{-2+o(1)} \to 0$  as $N\rightarrow \infty$.
Summing up, we have that 
\begin{equation}\label{casesmallm}
\sum_{m=2}^{2\ell_* +1}\bbP(\exists \text{ a component of size $m$}) \leq (2\ell_* +1) \varphi_2(N)\, +\, o(1)\rightarrow 0 \quad \text{as} \quad N\rightarrow \infty
\end{equation}
as desired.
 
\end{proof}

\smallskip
\subsection{Proof of Theorem \ref{connectedness} (1)--(3)}
Theorem \ref{connectedness} (1) is an immediate consequence of  Lemma \ref{zeroone-iso}. 

To prove Theorem \ref{connectedness} (2)--(3), we observe first the following fact:  Let $C$ be the event that $G(N,\K^r,p)$ is connected, $B_1$ the event that $G(N,\K^r,p)$ has isolated super-vertices  and $B_2$ the event that $G(N,\K^r,p)$ has components of size of between 2 and $N/2$. Hence,
\begin{eqnarray*}
\bbP(C)=\bbP(B_1^c)-\bbP(B_2)+\bbP(B_1\cap B_2).
\end{eqnarray*}

Let $X$ denote the number of isolated super-vertices in $G(N,\K^r,p)$.
Since $\bbP(B_1^c)=\bbP(X=0)$,  $\bbP(B_2)=o(1)$ because  Lemma \ref{smallcomponents}, and   $\bbP(B_1\cap B_2)<\bbP(B_2)=o(1)$. Thus, 
\begin{eqnarray*}
\bbP(C)=\bbP(G(N, K, p) \text{ is connected})&=&\bbP(X=0)+o(1),
%&\rightarrow& 1,
\end{eqnarray*}
%as $N\rightarrow\infty$ by Lemma \ref{zeroone-iso}.
Thus,  Theorem \ref{connectedness} (2) follows from Lemma \ref{poisson} and Lemma \ref{smallcomponents} while Theorem \ref{connectedness} (3) follows from Lemma  \ref{zeroone-iso} and Lemma \ref{smallcomponents}. 

% and Theorem \ref{connectedness} (3) from Lemma  \ref{zeroone-iso} and Lemma \ref{smallcomponents}. 

%%%%%%%%%%%%%%%%%%%%%%%%%%%%%%%%%%%%%%%
%%%%%%%%%%%%%%%%%%%%%%%%%%%%%%%%%%%%%%%%
\bigskip
\section{Proof of {\color{black} Proposition} \ref{corBollobas1} and {\color{black} Proposition} \ref{corBollobas}}\label{graphicalproof}

Consider the random graph $G(N,\K^r,p)$ with  $p=c/n$, where $c$ is a constant. Also assume  that the limit $u:=\lim_{N\to \infty}(n/N)$ exists. As discussed in Section \ref{relatedwork}, the random graph model $G(N,\K^r,p)$ belongs to the class of IRG studied in \cite{BJR},  where the emergence of the giant component and the degree distribution were analyzed under some  conditions on the sequence of kernels 
$\{\kappa_N\}_{N\geq1}$ given by \eqref{kernel}. Let us begin with some definitions. 

\bigskip
\begin{defn}[Definition 2.9 \cite{BJR}]
Consider a set of $N$ vertices.  Let $\mathcal{V}=({\mathcal S},\mu,({\bf x}_N)_{N\geq 1})$ be a vertex space and let $\kappa$ be a kernel on ${\mathcal V}$. 
A sequence $\{\kappa_N\}_{N\geq1}$ of kernels is graphical on $({\mathcal S},\mu)$ with limit $\kappa$ if the following holds.
\begin{enumerate}
\item[(1)] \label{condconv} For a.e. $(x,y)\in {\mathcal S}^2$, $x_N \rightarrow x  \text{ and }y_N \rightarrow y$
imply that  $\kappa_N(x_N,y_N)\rightarrow \kappa(x,y)$.
\item[(2)] \label{condcont} $\kappa$ is continuous a.e. on ${\mathcal S}^2$.
\item[(3)] \label{condL1}$\kappa \in L^1({\mathcal S}^2,\mu\times\mu)$.
\item[(4)] \label{condelos}If $e(G)$ is the number of edges of $G^{\mathcal V}(N,\kappa_N)$, then
$$\lim_{N\to \infty}\frac{\bbE(e(G))}{N}= \frac{1}{2}\int\int_{{\mathcal S}^2}\kappa(x,y)d\mu(x) d\mu(y).$$
\end{enumerate}
\end{defn}

\bigskip
Note that whether $\{\kappa_N\}_{N\geq1}$ is graphical depends on the sequences $\{{\bf x}_N\}_{N\geq 1}$.  The next lemma says that the only condition in our random graph model $G(N,\K^r,p)$ for $\{\kappa_N\}_{N\geq1}$ to be graphical is that the limit $u$ exists.

\bigskip
\begin{lem}\label{graphical}
Consider the random graph $G(N,\K^r,p)$ with $p=c/n$, where $c$ is a constant. If the limit $u$ exists, then the sequence of kernels $\{\kappa_N\}_{N\geq1}$ given by \eqref{kernel}  is graphical on $\mathcal V$ with limit given by the kernel $\kappa$ such that $\kappa(i,j)=(c/u)ij$, for $i,j\in \mathcal{S}$.
\end{lem}

\begin{proof}
It is not difficult to see that $\kappa$ is in fact a kernel on ${\mathcal V}$, and that the conditions (2) and (3) in Definition 1 are satisfied. To check the condition $(1)$ in Definition 1, we fix a point $(i,j)$ of $\mathbb{Z}^{+} \times \mathbb{Z}^{+}$ and consider two sequences, $(i_N)_{N\geq 1}$ and $(j_N)_{N\geq 1}$, such that $i_N \rightarrow i$ and $j_N \rightarrow j$. Since
$$\kappa_N(i_N,j_N)= \frac{N}{n}c \,i_N \,j_N + o(1),$$
and $u=\lim_{N\to \infty}(n/N)$ exists, we conclude that
$$\lim_{N\to \infty}\kappa(i_N,j_N)= \left(\frac{c}{u}\right)ij.$$
To show that the condition (4) is satisfied, note first that in $G(N,\K^r,p)$, 
$$ (1/2)\int\int_{{\mathcal S}^2}\kappa(x,y)d\mu(x) d\mu(y) = \frac{1}{2}\sum_{j=1}^{r}\sum_{i=1}^{r}\kappa(i,j)\mu_i \mu_j = \frac{c}{2u} \sum_{j=1}^{r} j \mu_j\sum_{i=1}^{r}i \mu_i = \frac{c\, u}{2}.$$
Furthermore, we have
\begin{eqnarray*}
\frac{\bbE(e(G))}{N}=\frac{1}{N}\sum_{1\leq k<l\leq N}p_{x_kx_l}&=&\frac{1}{N}\sum_{1\leq k<l\leq N}\frac{\kappa_N(x_k,x_l)}{N}\\
&=&\frac{1}{N^2}\sum_{1\leq k<l\leq N}\frac{N}{n}c\,x_k \,x_l + o(1)\\
&\leq&\frac{N}{n}\frac{c}{2}\left(\frac{1}{N}\sum_{k=1}^N x_k\right)^2 + o(1)\\
&=&\frac{N}{n}\frac{c}{2}\left(\frac{n}{N}\right)^2 + o(1).
\end{eqnarray*}
Therefore, we get
$$
\lim_{N\to\infty} \frac{\bbE(e(G))}{N}\leq \frac{c\,u}{2}.
$$
Finally, by Lemma 8.1 (\cite{BJR}) we know that if $\kappa$ is a continuous kernel on a vertex space ${\mathcal V}$, then
\begin{eqnarray*}
\liminf_{N\rightarrow\infty}\frac{\bbE(e(G))}{N}\geq \frac{1}{2}\sum_{j=1}^{\infty}\sum_{i=1}^{\infty}\kappa(i,j)\mu_i \mu_j.
\end{eqnarray*}
Thus, we have  
\begin{eqnarray*}
\lim_{N\to \infty}\frac{\bbE (e(G))}{N} = \frac{1}{2}\sum_{j=1}^{\infty}\sum_{i=1}^{\infty}\kappa(i,j)\mu_i \mu_j
\end{eqnarray*}
and the sequence of kernels $\{\kappa_N\}_{N\geq 1}$ is graphical on ${\mathcal V}$. 
\end{proof}

\bigskip
In order to obtain results concerning the size of the giant component, one additional definition is required.

\bigskip
\begin{defn}[Definition 2.10-11 \cite{BJR}]\label{irreducible}
A kernel $\kappa$ on $({\mathcal S},\mu)$ is irreducible if for all $ A\subset {\mathcal S}$ and $\kappa=0$ a.e. on $A\times(S\setminus A)$ implies $\mu(A)=0$ or $\mu(S\setminus A)=0$. 
\end{defn}

\bigskip
In fact, for technical reasons, a slight weakening of irreducibility is considered.

 \bigskip
\begin{defn}[Definition 2.11 \cite{BJR}]\label{quasiirreducible}
A kernel $\kappa$ on $({\mathcal S},\mu)$ is quasi-irreducible if there is a $\mu$-continuity set ${\mathcal S}'\subseteq{\mathcal S}$ with $\mu({\mathcal S}')>0$ such that the restriction of $\kappa$ to ${\mathcal S}'\times{\mathcal S}'$ is irreducible, and $\kappa(x,y)=0$ if $x\notin{\mathcal S}'$ or $y\notin{\mathcal S}'$. 
\end{defn}

\bigskip
Now, Propositions \ref{corBollobas1} and \ref{corBollobas} can be obtained as corollaries of Theorem 3.1 and Theorem 3.13 in (\cite{BJR}), respectively, which are included here for the sake of completeness.

\bigskip
\begin{teo}\label{T3.1}(Theorem 3.1 \cite{BJR})
Let $\{\kappa_N\}_{N\geq1}$ be a graphical sequence of kernels on a  vertex space ${\mathcal V}$ with limit $\kappa$, and let $G_N:=G^{\mathcal V}(N,\kappa_N)$.   Let $T_{\kappa}$ be an integral operator defined by
$$(T_{\kappa}f)(x)=\int_{{\mathcal S}}\kappa(x,y)f(y)d\mu(y),$$
for any (measurable) function $f$, such that its integral is defined (finite or $+\infty$) for a.e. $x$. Let $||T_{\kappa}||:=\sup\{||T_{\kappa}f||_{2}:f\geq 0 \text{ and }||f||_{2}\leq 1\}.$
\begin{enumerate}
\item[(1)] If $||T_{\kappa}||\leq1$, then 
$$\lim_{N\to \infty}\frac{L_1(G_N)}{N}=0$$
in probability, while if $||T_{\kappa}||>1$, then $whp$ $L_1(G_N)=\Theta(N)$.
\item[(2)] For any $\epsilon>0$, $whp$
$$\frac{L_1(G_N)}{N}\leq \rho(\kappa)+\epsilon.$$
\item[(3)] If $\kappa$ is quasi-irreducible, then
$$\lim_{N\to \infty}\frac{L_1(G_N)}{N} = \rho(\kappa),$$
in probability, where $\rho(\kappa):=\int_{\mathcal S}\rho(\kappa,x)d\mu(x),$ and the function $x \mapsto \rho(\kappa,x)$ is the maximal fixed point of the non-linear operator $\Phi_{\kappa}$ defined by $\Phi_{\kappa} f:=1-e^{-T_{\kappa}f}$. In addition, $\rho(\kappa)<1$, and $\rho(\kappa)>0$ if and only if $||T_{\kappa}||>1$. 
\end{enumerate}
\end{teo}

More details regarding the definition and behavior of $\rho(\kappa)$ can be found in Theorem 6.2, \cite{BJR}.

\bigskip
\begin{teo}\label{T3.13}(Theorem 3.13 \cite{BJR})
Let $\{\kappa_N\}_{N\geq1}$ be a graphical sequence of kernels on a  vertex space ${\mathcal V}$ with limit $\kappa$, and let $G_N:=G^{\mathcal V}(N,\kappa_N)$.  Let $Z_k$ the number of vertices of $G_N$ with degree $k$, for $k\geq 0$.  Then, for any fixed $k \geq 0$, 
$$\lim_{N\to \infty}\frac{Z_k}{N} =\bbP(\Xi= k),$$
in probability, where $\Xi$ has the mixed Poisson distribution $\int_{{\mathcal S}} Po(\lambda(x))d\mu(x)$, and $\lambda(x):=\int_{{\mathcal S}}\kappa(x,y) d\mu(y)$.
\end{teo}

\smallskip
\subsection{Proof of {\color{black} Proposition} \ref{corBollobas1}}

Observe that the kernel $\kappa$ has the form $\kappa(i,j)=\varphi(i)\varphi(j)$, with $\varphi(i):= (c/u)^{1/2} i$, which is the rank-$1$ case studied in \cite{BJR}. In this case, we have
$$||T_{\kappa}||=\int_{{\mathcal{S}}}\varphi^2 d\mu=\sum_{i=1}^{r} \frac{c}{u} i^2 \mu_i=c\bar{s_2}.$$
On the other hand, note that $\kappa(i,j)=\left(c/u\right) ij$ defined  on $(\bbZ^+,\mu)$, where $\mu(\{i\})=\mu_i$ given by \eqref{convki}, equals 0 only if $i=0$ or $j=0$. Hence,  $\kappa$ is irreducible and  quasi-irreducible.  
Therefore, {\color{black} Proposition} \ref{corBollobas1} follows as a consequence of Theorem \ref{T3.1} and Lemma \ref{graphical}.

\smallskip
\subsection{Proof of {\color{black} Proposition} \ref{corBollobas}}

{\color{black} Proposition} \ref{corBollobas} follows as a consequence of Theorem \ref{T3.13} and  Lemma \ref{graphical}.

\bigskip
\section{Discussions}\label{sec:diss}

In this section we compare our  results with related work.
 
\subsection{Comparison with connectedness of $G(n,p)$} In order to compare Theorem \ref{connectedness} with the threshold for connectedness of $G(n,p)$, take $p=(1/n)(\ln n+c(n))$. It is well known (see for example \cite{Bollobas}) that if $\lim_{n\rightarrow\infty}c(n)=-\infty$, then $whp$ $G(n,p)$ is disconnected, but if $\lim_{n\rightarrow\infty}c(n)=+\infty$, then $whp$ $G(n,p)$ is connected. Furthermore,  if $\lim_{n\rightarrow\infty}c(n)= c$ is a constant,  
\begin{eqnarray*} 
 \lim_{n\rightarrow\infty} \bbP[\,G(n,p) \text{ is connected }]=e^{-e^{-c}}.
\end{eqnarray*}
 
In  Theorem \ref{connectedness} (2) we assume that $\lim_{N\rightarrow\infty}(n/N)=1$. Since $N:=\sum_{i=1}^rk_{i}$ and $n:=\sum_{i=1}^rik_{i}$,  we have
$$\frac{n}{N}=1+\frac{1}{N}\sum_{i=2}^r(i-1)k_{i}.$$ Because $\lim_{N\rightarrow\infty}n/N=1$, $\lim_{N\rightarrow\infty}k_i/N=0$ for each $i=2,3,\ldots,r$. By \eqref{convki} we have $\mu_1=1$ and  $\mu_i=0$ for each $i\in \{2,3,\ldots,r\}$, which implies  $\ell_*=1$, where $\ell_*:=\min\{\ell \in\{1,2,\ldots,r\} : \mu_{\ell}>0\}$.
By the definition of $i_*:=\Big\{ 1\le i \le r \ \Big|\  k_{i}^{\frac{1}{i (n-i)}} = \max_{1\le j\le r} k_{j}^{\frac{1}{j (n-j)}}\Big\}$, we have $k_{i_*}\geq k_{1}^{\frac{i_* (n-i_*)}{n-1}} = k_{1}^{\frac{i_*}{(1+o(1))}}$. On the other hand, $\mu_1=1$ means $k_1= (1+o(1)) N$. Thus  $i_*$ should coincide with $\ell_*=1$. Furthermore, we have $\gamma(1)= \lim_{N\to \infty} \sum_{i\neq 1} k_i\, (k_{1}e^{c})^{-\frac{i(n-i)}{n-1}} =0$. 

Therefore Theorem \ref{connectedness} (2) yields  
\begin{eqnarray*} 
 \lim_{N\rightarrow\infty} \bbP[\,G(N,\K^r,p) \text{ is connected }]=e^{-e^{-c}}.
\end{eqnarray*}
In words, the asymptotic probability of $G(N,\K^r,p)$ being connected is the same as that of $G(n,p)$ being connected, provided $\lim_{N\rightarrow\infty}(n/N)=1$.
This is not surprising because in this case both models have the same asymptotic number of vertices.\\

\smallskip
\subsection{Comparison with connectedness of IRG}\label{sec:devroye} The property of inhomogeneous random graphs being connected has recently been studied by Devroye and Fraiman in \cite{Devroye}, who considered a random graph model on a set of $N$ vertices described in a similar way as the general model of \cite{BJR}. In such a model, each pair of vertices, say $k$ and $l$, are connected independently with probability
\begin{equation*}
p_{kl}:=\min\{1,\kappa(x_k,x_l)\}p_N,
\end{equation*}
where $\kappa$ is a kernel on the respective vertex space $\mathcal{V}$ (see Section \ref{relatedwork}),  $x_k,x_l$ are values associated with the vertices $k$ and $l$, respectively and $p_N := (\log\, N)/ N$.  In \cite{Devroye}, a connectivity threshold is obtained in terms of an \textit{isolation parameter}  $\lambda_*:=\text{ess} \inf \lambda (x)$, where
$$\lambda(x):=\int_{\mathcal{S}}\kappa(x,y)d\mu(y),$$
and $\text{ess} \inf \lambda (x) := \sup\{a\in \mathbb{R}:\mu\{x:\lambda(x)<a\}=0\}$. More precisely, it is proved that when $\lambda_*>1$ the graph is connected $whp$, while when $\lambda_*<1$ the graph is disconnected $whp$. 

In order to compare  our model $G(N,\K^r,p)$ to that of Devroye and Fraiman \cite{Devroye}, we consider $G(N,\K^r,p)$ with $r\in \bbN$ constant and $p:=\frac{\ln k_{i_{*}}+c(N)}{i_{*}(n-\ell_*)}$ such that $i_*=\ell_*$, i.e,  $\mu_{i_*}>0$. Thus  $k_{i_*}=\mu_{i_*} N +o(N)$ and so without loss of generality we may $c(N) = o(\ln N)$. By \eqref{prob}, we have, for $i,j=1,\ldots, r$,
\begin{align*}
p_{ij}&:= 1-\left(1-p\right)^{ij}=i\,j\,p + O(p^2)\\
& = \left(\frac{i\,j}{\, i_*}\right)\left(\frac{\ln N}{n-i_*}\right) + o(1) + O(p^2) \sim \left(\frac{i\,j}{u\, i_*}\right)\frac{\ln N}{N},
\end{align*}
where $u:=\lim_{N\to \infty} (n/N)$. So, if we consider the kernel 
$$\kappa(i,j)=\frac{i\,j}{u\, i_*},$$ then the connection probability between pairs of super-vertices in our model is approximately the same as the connection probability between pairs of vertices in the IRG introduced by \cite{Devroye}, provided $N$ is sufficiently large. Moreover, this choice of kernel corresponds to the case $\lambda_*= 1$, because  
\begin{align*}
\lambda(i)=\sum_{j=1}^{r}\kappa(i,j)\mu_j  
&=\sum_{j=1}^{r}\frac{i\,j}{u\, i_*}  \lim_{N\to \infty}\frac{k_{j}}{N}\\
&= \frac{i}{i_*}      \sum_{j=1}^{r} \frac{j }{\lim_{N\to \infty} (n/N)} \lim_{N\to \infty}\frac{k_{j}}{N}\\
& =  \frac{i}{i_*} \sum_{j=1}^{r} \frac{jk_j}{n}   = \frac{i}{i_*}
\end{align*} 
and $\sum_{i<a i_*}\mu_i = 0$ if and only if $a<1$. 
The case $\lambda_* =1$ was not covered by \cite{Devroye}.

\smallskip
\subsection{Unbounded sizes}  The proof of Theorem \ref{connectedness} relies mainly on the analysis of the asymptotic  number of isolated super-vertices in $G(N,\K^r,p)$. 
Note that if the number $r$ of sizes of super-vertices is either a constant independent of $N$ or $r=r(N)$ is tending to a constant as $N\to \infty$, $\ell_*:=\min\{i\in\{1,\ldots,r\}:\mu_i>0\}$,   (defined as  in the proof of Lemma 7 and Lemma 8), is well defined. Moreover   $\lim_{N\to \infty} (n/N)$  exists. However, if $r=r(N)$ satisfies $r(N)\rightarrow \infty$ as $N \rightarrow \infty$, we can not guarantee neither the existence of an integer $i\in \{1,2,3,\ldots\}$ with $\mu_i>0$ nor that of the  limit $\lim_{N\to \infty} (n/N)$. In this case the exact distribution of sizes of super-vertices is related to the existence of isolated super-vertices in a more complex manner.

In Theorem \ref{connectedness} we considered $r$  either a constant independent of $N$ or $r=r(N)$ tending to a constant as $N\to \infty$. However our analysis can be extended to  the case that  $r=r(N) \rightarrow \infty$ as $N \rightarrow \infty$.   The proof may follow the lines of the proof of  Theorem \ref{connectedness} with a few  modifications and possibly other additional conditions.

%e.g. in Theorem \ref{connectedness} (2) the existence of the limit
%$$\lim_{N\to \infty} \bbE[X] = \lim_{N\to \infty} \sum_{i=1}^{r(N)} k_i\, (k_{i_*}e^{c})^{-\frac{i(n-i)}{i_{*}(n-i_*)}}.$$

%\bigskip
\section*{Acknowledgements}
We thank Serguei Popov for suggesting us the construction of the inhomogeneous random graph model studied in the paper. We also thank Luiz Renato Fontes for fruitful discussions during the early stages of this work. The first two authors were financially supported by DFG KA 2748/3-1 and the Austrian Science Fund (FWF): P26826, and the last one by FAPESP 2013/03898-8, 2015/03868-7 and CNPq 479313/2012-1. The second and the third author also thank, respectively, ICMC - Universidade de São Paulo and Università di Torino, for their hospitality. Finally, we thank the referees for their careful reading of the manuscript and many valuable comments and suggestions that have helped to improve the paper.

%%%%%%%%%%%%%%%%%%%%%%%%%%%%%%%%%%%%%%%%%%%%%%%%%%%%%%%%%%%%%%%%%%%%%%%%%%%%%%
%%%%%%%%%%%%% REFERENCES
%%%%%%%%%%%%%%%%%%%%%%%%%%%%%%%%%%%%%%%%%%%%%%%%%%%%%%%%%%%%%%%%%%%%%%%%%%%%%%

\end{document}